\documentclass[a4paper,12pt,fleqn]{amsart}%
\usepackage{amsfonts}
\usepackage{amsmath}
\usepackage{amssymb}
\usepackage[pdftex]{graphicx}
\usepackage[a4paper,left=2cm,right=1.9cm,top=2cm,bottom=2cm]{geometry}%
\providecommand{\U}[1]{\protect\rule{.1in}{.1in}}
\newtheorem{theorem}{Theorem}
\theoremstyle{plain}

\newtheorem{definition}{Definition}

\newtheorem{lemma}{Lemma}

\newtheorem{proposition}{Proposition}
\newtheorem{remark}{Remark}

\numberwithin{equation}{section}
\begin{document}
\title[ ]{Some Results on Simpson type Conformable fractional inequalities}
\author[Z. \c{S}anl\i\ ]{Zeynep \c{S}anl\i}
\address{Department of Mathematics, Faculty of Sciences, Karadeniz Technical
University, 61080, Trabzon, Turkey}
\email{zeynep.sanli@ktu.edu.tr}
\thanks{$^{\ast}$Corresponding author}

\begin{abstract}
In this paper we established a new Simpson type conformable fractional
integral equality for convex functions. Based on this identity, some results
related to Simpson-like type inequalities are obtained. These results are then
applied to some special means of real numbers and two special functions,
modified Bessel function and $q-$digamma function, respectively.

\end{abstract}
\maketitle

\section{Introduction}

We will start with the following inequality is well known in the literature as
Simpson's inequality.

\begin{theorem}
\label{1.1}Let $\psi:\left[  \gamma,\delta\right]  \rightarrow%
%TCIMACRO{\U{211d} }%
%BeginExpansion
\mathbb{R}
%EndExpansion
$ be a four times continuously differentiable mapping on $\left(
\gamma,\delta\right)  $ and $\left\Vert \psi^{\left(  4\right)  }\right\Vert
_{\infty}=\sup\left\vert \psi^{\left(  4\right)  }\left(  \varepsilon\right)
\right\vert <\infty.$ Then, the following inequality holds:%
\begin{equation}
\left\vert
%TCIMACRO{\dint \limits_{\gamma}^{\delta}}%
%BeginExpansion
{\displaystyle\int\limits_{\gamma}^{\delta}}
%EndExpansion
\psi\left(  \varepsilon\right)  d\varepsilon-\frac{\delta-\gamma}{3}\left[
\frac{\psi\left(  \gamma\right)  +\psi\left(  \delta\right)  }{2}+2\psi\left(
\frac{\gamma+\delta}{2}\right)  \right]  \right\vert \leq\frac{1}%
{2880}\left\Vert \psi^{\left(  4\right)  }\right\Vert _{\infty}.\left(
\delta-\gamma\right)  ^{4}. \label{1-1}%
\end{equation}

\end{theorem}

This inequality \ref{1-1} has been studied by several authors, these papers
can be seen in
\cite{ADD09,DAC00,LDKZ18,M15,M17,RAJNN19,SB19,SSO10,SSO10a,TGB18,ZWD20}.

In \cite{SSO10}, Sarikaya et al. obtained the following inequality for
differentiable convex functions on Simpson's inequality and used the following
lemma to show this.

\begin{lemma}
\label{1.2}Let $\psi:I\subset%
%TCIMACRO{\U{211d} }%
%BeginExpansion
\mathbb{R}
%EndExpansion
\rightarrow%
%TCIMACRO{\U{211d} }%
%BeginExpansion
\mathbb{R}
%EndExpansion
$ be an absolutely continous mapping on $I^{\circ}$ such that $\psi^{\prime
}\in L_{1}\left[  \gamma,\delta\right]  ,$ where $\gamma,\delta\in I^{\circ}$
\ with $\gamma<\delta.$ Then, the following equality holds:%
\begin{align}
&  \frac{1}{6}\left[  \psi\left(  \gamma\right)  +4\psi\left(  \frac
{\gamma+\delta}{2}\right)  +\psi\left(  \delta\right)  \right]  -\frac
{1}{\delta-\gamma}%
%TCIMACRO{\dint \limits_{\gamma}^{\delta}}%
%BeginExpansion
{\displaystyle\int\limits_{\gamma}^{\delta}}
%EndExpansion
\psi\left(  \varepsilon\right)  d\varepsilon\label{1-2}\\
&  =\frac{\delta-\gamma}{2}%
%TCIMACRO{\dint \limits_{0}^{1}}%
%BeginExpansion
{\displaystyle\int\limits_{0}^{1}}
%EndExpansion
\left[  \left(  \frac{w}{2}-\frac{1}{3}\right)  \psi^{\prime}\left(
\frac{1+w}{2}\delta+\frac{1-w}{2}\gamma\right)  +\left(  \frac{1}{3}-\frac
{w}{2}\right)  \psi^{\prime}\left(  \frac{1+w}{2}\gamma+\frac{1-w}{2}%
\delta\right)  \right]  dw.\nonumber
\end{align}

\end{lemma}

The main theorem in \cite{SSO10} is as follows.

\begin{theorem}
\label{1.3}Let $\psi:I\subset%
%TCIMACRO{\U{211d} }%
%BeginExpansion
\mathbb{R}
%EndExpansion
\rightarrow%
%TCIMACRO{\U{211d} }%
%BeginExpansion
\mathbb{R}
%EndExpansion
$ be a differentiable mapping on $I^{\circ}$ such that $\psi^{\prime}\in
L_{1}\left[  a,b\right]  ,$ where $\gamma,\delta\in I^{\circ}$ \ with
$\gamma<\delta.$ If \ $\left\vert \psi^{\prime}\right\vert ^{q}$ is convex on
$\left[  \gamma,\delta\right]  ,q>1,$ then the following inequality holds:%
\begin{align}
&  \left\vert \frac{1}{6}\left[  \psi\left(  \gamma\right)  +4\psi\left(
\frac{\gamma+\delta}{2}\right)  +\psi\left(  \delta\right)  \right]  -\frac
{1}{\delta-\gamma}%
%TCIMACRO{\dint \limits_{\gamma}^{\delta}}%
%BeginExpansion
{\displaystyle\int\limits_{\gamma}^{\delta}}
%EndExpansion
\psi\left(  \varepsilon\right)  d\varepsilon\right\vert \label{1-3}\\
&  \leq\frac{\delta-\gamma}{12}\left(  \frac{1+2^{p+1}}{3\left(  p+1\right)
}\right)  ^{\frac{1}{p}}\left[  \left(  \frac{3\left\vert \psi^{\prime}\left(
\delta\right)  \right\vert ^{q}+\left\vert \psi^{\prime}\left(  \gamma\right)
\right\vert ^{q}}{4}\right)  ^{\frac{1}{q}}+\left(  \frac{\left\vert
\psi^{\prime}\left(  \delta\right)  \right\vert ^{q}+3\left\vert \psi^{\prime
}\left(  \gamma\right)  \right\vert ^{q}}{4}\right)  ^{\frac{1}{q}}\right]
.\nonumber
\end{align}
where $\frac{1}{p}+\frac{1}{q}=1.$
\end{theorem}

\begin{definition}
\label{1.4}A function $\psi:\left[  \gamma,\delta\right]  \rightarrow%
%TCIMACRO{\U{211d} }%
%BeginExpansion
\mathbb{R}
%EndExpansion
$ is said to be convex on $\left[  \gamma,\delta\right]  $ if the inequality%
\[
\psi\left(  wa+\left(  1-w\right)  b\right)  \leq wf\left(  a\right)  +\left(
1-w\right)  f\left(  b\right)
\]
holds for all $a,b\in\left[  \gamma,\delta\right]  $ and $w\in\left[
0,1\right]  .$ If $\left(  -\psi\right)  $ is convex, $\psi$ is concave.
\end{definition}

Convex functions are important for mathematical inequalities. Many authors
obtained several inequalities for convex functions
\cite{ASKI16,KIAS18,SK17,SNT15}. The most famous inequality has been used with
convex functions is Hermite-Hadamard, which is stated as follows:

Let $\psi:\left[  \gamma,\delta\right]  \rightarrow%
%TCIMACRO{\U{211d} }%
%BeginExpansion
\mathbb{R}
%EndExpansion
$ be a convex function and $a,b\in\left[  \gamma,\delta\right]  $ with $a<b.$
Then the following double inequalities hold:%
\[
\psi\left(  \frac{a+b}{2}\right)  \leq\frac{1}{b-a}%
%TCIMACRO{\dint \limits_{a}^{b}}%
%BeginExpansion
{\displaystyle\int\limits_{a}^{b}}
%EndExpansion
\psi\left(  \varepsilon\right)  d\varepsilon\leq\frac{\psi\left(  a\right)
+\psi\left(  b\right)  }{2}.
\]

The aim of this paper is to establish Simpson type conformable fractional
integral inequalities based on convexity.

\section{Preliminaries}

In this section, we give some definitions and basic results we will use.

\begin{definition}
\label{2.1} Let $\gamma,\delta\in%
%TCIMACRO{\U{211d} }%
%BeginExpansion
\mathbb{R}
%EndExpansion
$ with $\gamma<\delta$ and $f\in L\left[  \gamma,\delta\right]  .$ The left
and right Riemann- Liouville fractional integrals $J_{\gamma+}^{\tau}f$ and
$J_{\delta-}^{\tau}f$ of order $\tau>0$ are defined by
\[
J_{\gamma+}^{\tau}f=\frac{1}{\Gamma\left(  \tau\right)  }%
%TCIMACRO{\dint \limits_{\gamma}^{\varepsilon}}%
%BeginExpansion
{\displaystyle\int\limits_{\gamma}^{\varepsilon}}
%EndExpansion
\left(  \varepsilon-w\right)  ^{\tau-1}\psi\left(  w\right)  dw,\text{
\ }\varepsilon>\gamma
\]
and%
\[
J_{\delta-}^{\tau}f=\frac{1}{\Gamma\left(  \tau\right)  }%
%TCIMACRO{\dint \limits_{\varepsilon}^{\delta}}%
%BeginExpansion
{\displaystyle\int\limits_{\varepsilon}^{\delta}}
%EndExpansion
\left(  w-\varepsilon\right)  ^{\tau-1}\psi\left(  w\right)  dw,\text{
\ }\varepsilon<\delta
\]
respectively, where $\Gamma\left(  \tau\right)  $ is the Gamma function
defined by $\Gamma\left(  \tau\right)  =%
%TCIMACRO{\dint \limits_{0}^{\infty}}%
%BeginExpansion
{\displaystyle\int\limits_{0}^{\infty}}
%EndExpansion
e^{-w}t^{\tau-1}dw$ (see \cite{KST06}, p. 69).
\end{definition}

The following defiinition of conformable fractional integrals could be found
in \cite{A15,SK17}.

\begin{definition}
\label{2.2} Let $\tau\in(m,m+1],$ $m=0,1,2,...,$ $\beta=\tau-m,$
$\gamma,\delta\in%
%TCIMACRO{\U{211d} }%
%BeginExpansion
\mathbb{R}
%EndExpansion
$ with $\gamma<\delta$ and $\psi\in L\left[  \gamma,\delta\right]  .$ The left
and right conformable fractional integrals $%
\begin{array}
[c]{c}%
I_{\tau}^{\gamma}\psi
\end{array}
$ and $%
\begin{array}
[c]{c}%
^{\delta}I_{\tau}\psi
\end{array}
$ of order $\alpha>0$ are defined by%
\[
I_{\tau}^{\gamma}\psi=\frac{1}{m!}%
%TCIMACRO{\dint \limits_{\gamma}^{\varepsilon}}%
%BeginExpansion
{\displaystyle\int\limits_{\gamma}^{\varepsilon}}
%EndExpansion
\left(  \varepsilon-w\right)  ^{m}\left(  w-\gamma\right)  ^{\beta-1}%
\psi\left(  w\right)  dw,\ \varepsilon>\gamma
\]
and%
\[
^{\delta}I_{\tau}\psi=\frac{1}{m!}%
%TCIMACRO{\dint \limits_{\varepsilon}^{\delta}}%
%BeginExpansion
{\displaystyle\int\limits_{\varepsilon}^{\delta}}
%EndExpansion
\left(  w-\varepsilon\right)  ^{m}\left(  \delta-w\right)  ^{\beta-1}%
\psi\left(  w\right)  dw,\text{ \ }\varepsilon<\delta
\]
respectively.
\end{definition}

It is easily seen that if one takes $\tau=m+1$ in the Definition \ref{2.2}
(for the left and right conformable fractional integrals), one has the
Definition \ref{2.1} (the left and rightRiemann-Liouville fractional
integrals) for $\tau\in%
%TCIMACRO{\U{2115} }%
%BeginExpansion
\mathbb{N}
%EndExpansion
.$

\section{Main Results}

Let's begin start the following Lemma which helps us to obtain the main results:

\begin{lemma}
\label{3.1}Let $\psi:I\subset\left(  0,\infty\right)  \rightarrow%
%TCIMACRO{\U{211d} }%
%BeginExpansion
\mathbb{R}
%EndExpansion
,$ be a differentiable function on $I^{\circ},\gamma,\delta\in I^{\circ}$ and
$\gamma<\delta.$ If $\psi^{\prime}$ $\in L\left[  \gamma,\delta\right]  ,$
then the following equality holds:
\begin{align}
&  \frac{1}{6}\left[  \psi\left(  \gamma\right)  +4\psi\left(  \frac
{\gamma+\delta}{2}\right)  +\psi\left(  \delta\right)  \right]  -\frac
{2^{\tau-1}}{\left(  \delta-\gamma\right)  ^{\tau}}\frac{\Gamma\left(
\tau+1\right)  }{\Gamma\left(  \tau-m\right)  }\left[
\begin{array}
[c]{c}%
I_{\tau}^{\gamma}\psi\left(  \frac{\gamma+\delta}{2}\right)
\end{array}
+%
\begin{array}
[c]{c}%
^{\delta}I_{\tau}\psi\left(  \frac{\gamma+\delta}{2}\right)
\end{array}
\right]  \label{3-1}\\
&  =\frac{\delta-\gamma}{2.m!}\frac{\Gamma\left(  \tau+1\right)  }%
{\Gamma\left(  \tau-m\right)  }\left\{
\begin{array}
[c]{c}%
%TCIMACRO{\dint \limits_{0}^{1}}%
%BeginExpansion
{\displaystyle\int\limits_{0}^{1}}
%EndExpansion
\left(  \frac{1}{3}\beta\left(  m+1,\tau-m\right)  -\frac{1}{2}\beta
_{w}\left(  m+1,\alpha-n\right)  \right)  \\
\times%
\begin{array}
[c]{c}%
\psi^{\prime}\left(  \frac{1+w}{2}\gamma+\frac{1-w}{2}\delta\right)
\end{array}
dw\\
-%
%TCIMACRO{\dint \limits_{0}^{1}}%
%BeginExpansion
{\displaystyle\int\limits_{0}^{1}}
%EndExpansion
\left(  \frac{1}{3}\beta\left(  m+1,\tau-m\right)  -\frac{1}{2}\beta
_{w}\left(  m+1,\alpha-n\right)  \right)  \\
\times%
\begin{array}
[c]{c}%
\psi^{\prime}\left(  \frac{1-w}{2}\gamma+\frac{1+w}{2}\delta\right)
\end{array}
dw
\end{array}
\right\}  .\nonumber
\end{align}

\end{lemma}

\begin{proof}
We start by considering the following computations which follows from change
of variables and using the definition of the conformable fractional integrals.%

\begin{align*}
I_{1}  &  =\frac{1}{m!}%
%TCIMACRO{\dint \limits_{0}^{1}}%
%BeginExpansion
{\displaystyle\int\limits_{0}^{1}}
%EndExpansion
\left(  \frac{1}{3}\beta\left(  m+1,\tau-m\right)  -\frac{1}{2}\beta
_{w}\left(  m+1,\tau-m\right)  \right)
\begin{array}
[c]{c}%
\psi^{\prime}\left(  \frac{1+w}{2}\gamma+\frac{1-w}{2}\delta\right)
\end{array}
dw\\
&  =\frac{1}{m!}%
%TCIMACRO{\dint \limits_{0}^{1}}%
%BeginExpansion
{\displaystyle\int\limits_{0}^{1}}
%EndExpansion
\left(
\begin{array}
[c]{c}%
\frac{1}{3}\beta\left(  m+1,\tau-m\right)
\begin{array}
[c]{c}%
\psi^{\prime}\left(  \frac{1+w}{2}\gamma+\frac{1-w}{2}\delta\right)
\end{array}
dw\\
-\frac{1}{2}\beta_{w}\left(  m+1,\tau-m\right)
\begin{array}
[c]{c}%
\psi^{\prime}\left(  \frac{1+w}{2}\gamma+\frac{1-w}{2}\delta\right)
\end{array}
dw
\end{array}
\right) \\
&  =\frac{1}{3.m!}\frac{2}{\delta-\gamma}\beta\left(  m+1,\tau-m\right)
\begin{array}
[c]{c}%
\psi\left(  \frac{1+w}{2}\gamma+\frac{1-w}{2}\delta\right)
\end{array}
|_{0}^{1}\\
&  -\frac{1}{2.m!}\left(
%TCIMACRO{\dint \limits_{0}^{1}}%
%BeginExpansion
{\displaystyle\int\limits_{0}^{1}}
%EndExpansion
\left(
%TCIMACRO{\dint \limits_{0}^{w}}%
%BeginExpansion
{\displaystyle\int\limits_{0}^{w}}
%EndExpansion
\varepsilon^{n}\left(  1-\varepsilon\right)  ^{\tau-m-1}d\varepsilon\right)
\begin{array}
[c]{c}%
\psi^{\prime}\left(  \frac{1+w}{2}\gamma+\frac{1-w}{2}\delta\right)
\end{array}
dw\right) \\
&  =\frac{2}{\delta-\gamma}\frac{1}{3}\frac{\Gamma\left(  \tau-m\right)
}{\Gamma\left(  \tau+1\right)  }\left(  \psi\left(  \frac{\gamma+\delta}%
{2}\right)  -f\left(  \gamma\right)  \right) \\
&  -\frac{1}{2.m!}\frac{2}{\gamma-\delta}\left[
\begin{array}
[c]{c}%
-\left(
%TCIMACRO{\dint \limits_{0}^{w}}%
%BeginExpansion
{\displaystyle\int\limits_{0}^{w}}
%EndExpansion
\varepsilon^{n}\left(  1-\varepsilon\right)  ^{\tau-m-1}d\varepsilon\right)
\begin{array}
[c]{c}%
\psi\left(  \frac{1+w}{2}\gamma+\frac{1-w}{2}\delta\right)
\end{array}
|_{0}^{1}\\
+%
%TCIMACRO{\dint \limits_{0}^{1}}%
%BeginExpansion
{\displaystyle\int\limits_{0}^{1}}
%EndExpansion
w^{n}\left(  1-w\right)  ^{\tau-m-1}%
\begin{array}
[c]{c}%
\psi\left(  \frac{1+w}{2}\gamma+\frac{1-w}{2}\delta\right)
\end{array}
dw
\end{array}
\right] \\
&  =\frac{2}{\gamma-\delta}\frac{1}{3}\frac{\Gamma\left(  \tau-m\right)
}{\Gamma\left(  \tau+1\right)  }\left(  \psi\left(  \frac{\gamma+\delta}%
{2}\right)  -f\left(  \gamma\right)  \right) \\
&  +\frac{2}{\gamma-\delta}\frac{1}{2}\frac{\Gamma\left(  \tau-m\right)
}{\Gamma\left(  \tau+1\right)  }\psi\left(  \gamma\right)  -\frac{1}{2}\left(
\frac{2}{\gamma-\delta}\right)  ^{\tau+1}%
\begin{array}
[c]{c}%
I_{\tau}^{\gamma}\psi\left(  \frac{\gamma+\delta}{2}\right)
\end{array}
\\
&  =\frac{2}{\gamma-\delta}\frac{\Gamma\left(  \tau-m\right)  }{\Gamma\left(
\tau+1\right)  }\left(  \frac{1}{6}\psi\left(  \gamma\right)  +\frac{1}{3}%
\psi\left(  \frac{\gamma+\delta}{2}\right)  \right)  -\frac{1}{2}\left(
\frac{2}{\gamma-\delta}\right)  ^{\tau+1}%
\begin{array}
[c]{c}%
I_{\tau}^{\gamma}\psi\left(  \frac{\gamma+\delta}{2}\right)
\end{array}
\end{align*}

and similarly%

\begin{align*}
I_{2}  &  =\frac{1}{m!}%
%TCIMACRO{\dint \limits_{0}^{1}}%
%BeginExpansion
{\displaystyle\int\limits_{0}^{1}}
%EndExpansion
\left(  \frac{1}{3}\beta\left(  m+1,\tau-m\right)  -\frac{1}{2}\beta
_{w}\left(  m+1,\tau-m\right)  \right)
\begin{array}
[c]{c}%
\psi^{\prime}\left(  \frac{1-w}{2}\gamma+\frac{1+w}{2}\delta\right)
\end{array}
dw\\
&  =\frac{2}{\delta-\gamma}\frac{\Gamma\left(  \tau-m\right)  }{\Gamma\left(
\tau+1\right)  }\left(  -\frac{1}{6}\psi\left(  \delta\right)  -\frac{1}%
{3}\psi\left(  \frac{\gamma+\delta}{2}\right)  \right)  -\frac{1}{2}\left(
\frac{2}{\delta-\gamma}\right)  ^{\tau+1}%
\begin{array}
[c]{c}%
^{\delta}I_{\tau}\psi\left(  \frac{\gamma+\delta}{2}\right)
\end{array}
.
\end{align*}

Thus, we can \ write%
\begin{align*}
I_{1}-I_{2}  &  =\frac{1}{3\left(  \delta-\gamma\right)  }\frac{\Gamma\left(
\tau-m\right)  }{\Gamma\left(  \tau+1\right)  }\left[  \psi\left(
\gamma\right)  +\psi\left(  \delta\right)  +4\psi\left(  \frac{\gamma+\delta
}{2}\right)  \right] \\
&  -\frac{1}{2}\left(  \frac{2}{\delta-\gamma}\right)  ^{\tau+1}\left[
\begin{array}
[c]{c}%
I_{\tau}^{\gamma}\psi\left(  \frac{\gamma+\delta}{2}\right)
\end{array}
+%
\begin{array}
[c]{c}%
^{\delta}I_{\tau}\psi\left(  \frac{\gamma+\delta}{2}\right)
\end{array}
\right]  ,
\end{align*}

which is completes the proof.
\end{proof}

\begin{remark}
\label{3.2}If we take $\tau=m+1$ in Lemma \ref{3.1}, we have the following
equality%
\begin{align}
&  \frac{1}{6}\left[  \psi\left(  \gamma\right)  +4\psi\left(  \frac
{\gamma+\delta}{2}\right)  +\psi\left(  \delta\right)  \right]  -\frac
{2^{\tau-1}}{\left(  \delta-\gamma\right)  ^{\tau}}\Gamma\left(
\tau+1\right)  \left[
\begin{array}
[c]{c}%
J_{\gamma^{+}}^{\tau}\psi\left(  \frac{\gamma+\delta}{2}\right)
\end{array}
+%
\begin{array}
[c]{c}%
J_{\delta^{-}}^{\tau}\psi\left(  \frac{\gamma+\delta}{2}\right)
\end{array}
\right]  \label{3-2}\\
&  =\frac{\delta-\gamma}{2}\left[
\begin{array}
[c]{c}%
%TCIMACRO{\dint \limits_{0}^{1}}%
%BeginExpansion
{\displaystyle\int\limits_{0}^{1}}
%EndExpansion
\left(  \frac{1}{3}-\frac{w^{\tau}}{2}\right)
\begin{array}
[c]{c}%
\psi^{\prime}\left(  \frac{1+w}{2}\gamma+\frac{1-w}{2}\delta\right)
\end{array}
dw\\
-%
%TCIMACRO{\dint \limits_{0}^{1}}%
%BeginExpansion
{\displaystyle\int\limits_{0}^{1}}
%EndExpansion
\left(  \frac{1}{3}-\frac{w^{\tau}}{2}\right)
\begin{array}
[c]{c}%
\psi^{\prime}\left(  \frac{1+w}{2}\delta+\frac{1-w}{2}\gamma\right)
\end{array}
dw
\end{array}
\right]  \nonumber
\end{align}
which is proved by Matloka in (\cite{M15}, Lemma 5).
\end{remark}

\begin{remark}
\label{3.2a}If we take $\tau=1$ in Remark \ref{3.2}, we have the following
equality
\begin{align}
&  \frac{1}{6}\left[  \psi\left(  \gamma\right)  +4\psi\left(  \frac
{\gamma+\delta}{2}\right)  +\psi\left(  \delta\right)  \right]  -\frac
{1}{\delta-\gamma}%
%TCIMACRO{\dint \limits_{\gamma}^{\delta}}%
%BeginExpansion
{\displaystyle\int\limits_{\gamma}^{\delta}}
%EndExpansion
\psi\left(  \varepsilon\right)  d\varepsilon\label{3-2a}\\
&  =\frac{\delta-\gamma}{2}%
%TCIMACRO{\dint \limits_{0}^{1}}%
%BeginExpansion
{\displaystyle\int\limits_{0}^{1}}
%EndExpansion
\left(  \frac{1}{3}-\frac{w}{2}\right)  \left[
\begin{array}
[c]{c}%
\psi^{\prime}\left(  \frac{1+w}{2}\gamma+\frac{1-w}{2}\delta\right)
\end{array}
-%
\begin{array}
[c]{c}%
\psi^{\prime}\left(  \frac{1+w}{2}\delta+\frac{1-w}{2}\gamma\right)
\end{array}
\right]  dw\nonumber
\end{align}
in Theorem \ref{1.2}.
\end{remark}

\begin{theorem}
\label{3.3}Let$\psi:I\subset\left(  0,\infty\right)  \rightarrow%
%TCIMACRO{\U{211d} }%
%BeginExpansion
\mathbb{R}
%EndExpansion
,$ be a differentiable function on $I^{\circ},\gamma,\delta\in I^{\circ}$ and
$\gamma<\delta.$ If $\psi^{\prime}$ $\in L\left[  \gamma,\delta\right]  $ and
$\left\vert \psi^{\prime}\right\vert $ is a convex function on $\left[
\gamma,\delta\right]  $, then the following inequality holds:
\begin{align}
&  \left\vert \frac{1}{6}\left[  \psi\left(  \gamma\right)  +4\psi\left(
\frac{\gamma+\delta}{2}\right)  +\psi\left(  \delta\right)  \right]
-\frac{2^{\tau-1}}{\left(  \delta-\gamma\right)  ^{\tau}}\frac{\Gamma\left(
\tau+1\right)  }{\Gamma\left(  \tau-m\right)  }\left[
\begin{array}
[c]{c}%
I_{\tau}^{\gamma}\psi\left(  \frac{\gamma+\delta}{2}\right)
\end{array}
+%
\begin{array}
[c]{c}%
^{\delta}I_{\tau}\psi\left(  \frac{\gamma+\delta}{2}\right)
\end{array}
\right]  \right\vert \label{3-3}\\
&  \leq\frac{\delta-\gamma}{2.m!}\frac{\Gamma\left(  \tau-m\right)  }%
{\Gamma\left(  \tau+1\right)  }Z_{1}\left(  \tau,m\right)  \left(  \left\vert
\psi^{\prime}\left(  \gamma\right)  \right\vert +\left\vert \psi^{\prime
}\left(  \delta\right)  \right\vert \right) \nonumber
\end{align}
where%
\[
Z_{1}\left(  \tau,m\right)  =%
%TCIMACRO{\dint \limits_{0}^{1}}%
%BeginExpansion
{\displaystyle\int\limits_{0}^{1}}
%EndExpansion
\left\vert \frac{1}{3}\beta\left(  m+1,\tau-m\right)  -\frac{1}{2}\beta
_{w}\left(  m+1,\tau-m\right)  \right\vert dw
\]
with $m=0,1,2,...$ and $\tau\in\left(  m,m+1\right]  .$
\end{theorem}

\begin{proof}
From Lemma \ref{3.1} and $\left\vert \psi^{\prime}\right\vert $ is convex, we
have%
\begin{align*}
&  \left\vert \frac{1}{6}\left[  \psi\left(  \gamma\right)  +4\psi\left(
\frac{\gamma+\delta}{2}\right)  +\psi\left(  \delta\right)  \right]
-\frac{2^{\tau-1}}{\left(  \delta-\gamma\right)  ^{\tau}}\frac{\Gamma\left(
\tau+1\right)  }{\Gamma\left(  \tau-m\right)  }\left[
\begin{array}
[c]{c}%
I_{\tau}^{\gamma}\psi\left(  \frac{\gamma+\delta}{2}\right)
\end{array}
+%
\begin{array}
[c]{c}%
^{\delta}I_{\tau}\psi\left(  \frac{\gamma+\delta}{2}\right)
\end{array}
\right]  \right\vert \\
&  \leq\frac{b-a}{2.m!}\frac{\Gamma\left(  \tau+1\right)  }{\Gamma\left(
\tau-m\right)  }\left\{
\begin{array}
[c]{c}%
%TCIMACRO{\dint \limits_{0}^{1}}%
%BeginExpansion
{\displaystyle\int\limits_{0}^{1}}
%EndExpansion
\left\vert \left(  \frac{1}{3}\beta\left(  m+1,\tau-m\right)  -\frac{1}%
{2}\beta_{w}\left(  m+1,\tau-m\right)  \right)  \right\vert \\
\times\left(  \left\vert
\begin{array}
[c]{c}%
\psi^{\prime}\left(  \frac{1+w}{2}\gamma+\frac{1-w}{2}\delta\right)
\end{array}
\right\vert +\left\vert
\begin{array}
[c]{c}%
\psi^{\prime}\left(  \frac{1-w}{2}\gamma+\frac{1+w}{2}\delta\right)
\end{array}
\right\vert \right)  dw
\end{array}
\right\}  \\
&  \leq\frac{b-a}{2.m!}\frac{\Gamma\left(  \tau+1\right)  }{\Gamma\left(
\tau-m\right)  }\left\{
\begin{array}
[c]{c}%
%TCIMACRO{\dint \limits_{0}^{1}}%
%BeginExpansion
{\displaystyle\int\limits_{0}^{1}}
%EndExpansion
\left\vert \left(  \frac{1}{3}\beta\left(  m+1,\tau-m\right)  -\frac{1}%
{2}\beta_{w}\left(  m+1,\tau-m\right)  \right)  \right\vert \\
\times\left(
\begin{array}
[c]{c}%
\frac{1+w}{2}\left\vert
\begin{array}
[c]{c}%
\psi^{\prime}\left(  \gamma\right)
\end{array}
\right\vert +\frac{1-w}{2}\left\vert
\begin{array}
[c]{c}%
\psi^{\prime}\left(  \delta\right)
\end{array}
\right\vert \\
+\frac{1-w}{2}\left\vert
\begin{array}
[c]{c}%
\psi^{\prime}\left(  \gamma\right)
\end{array}
\right\vert +\frac{1+w}{2}\left\vert
\begin{array}
[c]{c}%
\psi^{\prime}\left(  \delta\right)
\end{array}
\right\vert
\end{array}
\right)  dw
\end{array}
\right\}  \\
&  =\frac{b-a}{2.m!}\frac{\Gamma\left(  \tau+1\right)  }{\Gamma\left(
\tau-m\right)  }%
%TCIMACRO{\dint \limits_{0}^{1}}%
%BeginExpansion
{\displaystyle\int\limits_{0}^{1}}
%EndExpansion
\left\vert \frac{1}{3}\beta\left(  m+1,\tau-m\right)  -\frac{1}{2}\beta
_{w}\left(  m+1,\tau-m\right)  \right\vert \\
&  \times\left(  \left\vert
\begin{array}
[c]{c}%
\psi^{\prime}\left(  \gamma\right)
\end{array}
\right\vert +\left\vert
\begin{array}
[c]{c}%
\psi^{\prime}\left(  \delta\right)
\end{array}
\right\vert \right)  dw\\
&  =\frac{b-a}{2.m!}\frac{\Gamma\left(  \tau+1\right)  }{\Gamma\left(
\tau-m\right)  }Z_{1}\left(  \tau,m\right)  \left(  \left\vert
\begin{array}
[c]{c}%
\psi^{\prime}\left(  \gamma\right)
\end{array}
\right\vert +\left\vert
\begin{array}
[c]{c}%
\psi^{\prime}\left(  \delta\right)
\end{array}
\right\vert \right)  .
\end{align*}

This completes the proof.
\end{proof}

\begin{remark}
\label{3.4} If we take $\tau=m+1,$ after that if we take $\tau=1$ in Theorem
\ref{3.3}, we obtain the following inequality
\end{remark}

\begin{align}
&  \left\vert \frac{1}{6}\left[  \psi\left(  \gamma\right)  +4\psi\left(
\frac{\gamma+\delta}{2}\right)  +\psi\left(  \delta\right)  \right]  -\frac
{1}{\delta-\gamma}%
%TCIMACRO{\dint \limits_{\gamma}^{\delta}}%
%BeginExpansion
{\displaystyle\int\limits_{\gamma}^{\delta}}
%EndExpansion
\psi\left(  \varepsilon\right)  d\varepsilon\right\vert \label{3-4}\\
&  \leq\frac{5\left(  b-a\right)  }{72}\left(  \left\vert f^{\prime}\left(
a\right)  \right\vert +\left\vert f^{\prime}\left(  b\right)  \right\vert
\right) \nonumber
\end{align}

where%

\[%
%TCIMACRO{\dint \limits_{0}^{1}}%
%BeginExpansion
{\displaystyle\int\limits_{0}^{1}}
%EndExpansion
\left\vert \frac{1}{2}w-\frac{1}{3}\right\vert dw=\frac{5}{36}.
\]
See also (\cite{ADD09}, Corollary 1).

\begin{theorem}
\label{3.6}Let $\psi:I\subset\left(  0,\infty\right)  \rightarrow%
%TCIMACRO{\U{211d} }%
%BeginExpansion
\mathbb{R}
%EndExpansion
,$ be a differentiable function on $I^{\circ},\gamma,\delta\in I^{\circ}$ and
$\gamma<\delta.$ If $\psi^{\prime}$ $\in L\left[  \gamma,\delta\right]  $ and
$\left\vert \psi^{\prime}\right\vert ^{q}$ is a convex function on $\left[
\gamma,\delta\right]  $ for $q>1$ and $\frac{1}{p}+\frac{1}{q}=1,$ then the
following inequality holds:
\begin{align}
&  \left\vert \frac{1}{6}\left[  \psi\left(  \gamma\right)  +4\psi\left(
\frac{\gamma+\delta}{2}\right)  +\psi\left(  \delta\right)  \right]
-\frac{2^{\tau-1}}{\left(  \delta-\gamma\right)  ^{\tau}}\frac{\Gamma\left(
\tau+1\right)  }{\Gamma\left(  \tau-m\right)  }\left[
\begin{array}
[c]{c}%
I_{\tau}^{\gamma}\psi\left(  \frac{\gamma+\delta}{2}\right)
\end{array}
+%
\begin{array}
[c]{c}%
^{\delta}I_{\tau}\psi\left(  \frac{\gamma+\delta}{2}\right)
\end{array}
\right]  \right\vert \label{3-6}\\
&  \leq\frac{\delta-\gamma}{m!}\frac{\Gamma\left(  \tau+1\right)  }%
{\Gamma\left(  \tau-m\right)  }\left(  \frac{1}{2}\right)  ^{2q+1}\left(
Z_{2}\left(  \tau,m\right)  \right)  ^{\frac{1}{p}}\left[
\begin{array}
[c]{c}%
\left(  3\left\vert
\begin{array}
[c]{c}%
\psi^{\prime}\left(  \gamma\right)
\end{array}
\right\vert ^{q}+\left\vert
\begin{array}
[c]{c}%
\psi^{\prime}\left(  \delta\right)
\end{array}
\right\vert ^{q}\right)  ^{\frac{1}{q}}\\
+\left(  \left\vert
\begin{array}
[c]{c}%
\psi^{\prime}\left(  \gamma\right)
\end{array}
\right\vert ^{q}+3\left\vert
\begin{array}
[c]{c}%
\psi^{\prime}\left(  \delta\right)
\end{array}
\right\vert ^{q}\right)  ^{\frac{1}{q}}%
\end{array}
\right]  .\nonumber
\end{align}
where%
\[
Z_{2}\left(  \tau,m\right)  =%
%TCIMACRO{\dint \limits_{0}^{1}}%
%BeginExpansion
{\displaystyle\int\limits_{0}^{1}}
%EndExpansion
\left\vert \frac{1}{3}\beta\left(  m+1,\tau-m\right)  -\frac{1}{2}\beta
_{w}\left(  m+1,\tau-m\right)  \right\vert ^{p}dw
\]
$m=0,1,2,...$ and $\tau\in\left(  m,m+1\right]  .$
\end{theorem}

\begin{proof}
From Lemma \ref{3.1} and using the H\"{o}lder's integral inequality and the
convexity of $\left\vert \psi^{\prime}\right\vert ^{q}$, we have
\begin{align*}
&  \left\vert \frac{1}{6}\left[  \psi\left(  \gamma\right)  +4\psi\left(
\frac{\gamma+\delta}{2}\right)  +\psi\left(  \delta\right)  \right]
-\frac{2^{\tau-1}}{\left(  \delta-\gamma\right)  ^{\tau}}\frac{\Gamma\left(
\tau+1\right)  }{\Gamma\left(  \tau-m\right)  }\left[
\begin{array}
[c]{c}%
I_{\tau}^{\gamma}\psi\left(  \frac{\gamma+\delta}{2}\right)
\end{array}
+%
\begin{array}
[c]{c}%
^{\delta}I_{\tau}\psi\left(  \frac{\gamma+\delta}{2}\right)
\end{array}
\right]  \right\vert \\
&  \leq\frac{\delta-\gamma}{2.m!}\frac{\Gamma\left(  \tau+1\right)  }%
{\Gamma\left(  \tau-m\right)  }\left\{
\begin{array}
[c]{c}%
\left(
%TCIMACRO{\dint \limits_{0}^{1}}%
%BeginExpansion
{\displaystyle\int\limits_{0}^{1}}
%EndExpansion
\left\vert \frac{1}{3}\beta\left(  m+1,\tau-m\right)  -\frac{1}{2}\beta
_{w}\left(  m+1,\tau-m\right)  \right\vert ^{p}dw\right)  ^{\frac{1}{p}}\\
\times\left(
%TCIMACRO{\dint \limits_{0}^{1}}%
%BeginExpansion
{\displaystyle\int\limits_{0}^{1}}
%EndExpansion
\left\vert
\begin{array}
[c]{c}%
\psi^{\prime}\left(  \frac{1+w}{2}\gamma+\frac{1-w}{2}\delta\right)
\end{array}
\right\vert ^{q}dw\right)  ^{\frac{1}{q}}\\
+\left(
%TCIMACRO{\dint \limits_{0}^{1}}%
%BeginExpansion
{\displaystyle\int\limits_{0}^{1}}
%EndExpansion
\left\vert \frac{1}{3}\beta\left(  m+1,\tau-m\right)  -\frac{1}{2}\beta
_{w}\left(  m+1,\tau-m\right)  \right\vert ^{p}dw\right)  ^{\frac{1}{p}}\\
\times\left(
%TCIMACRO{\dint \limits_{0}^{1}}%
%BeginExpansion
{\displaystyle\int\limits_{0}^{1}}
%EndExpansion
\left\vert
\begin{array}
[c]{c}%
\psi^{\prime}\left(  \frac{1-w}{2}\gamma+\frac{1+w}{2}\delta\right)
\end{array}
\right\vert ^{q}dw\right)  ^{\frac{1}{q}}%
\end{array}
\right\}  \\
&  \leq\frac{\delta-\gamma}{2.m!}\frac{\Gamma\left(  \tau+1\right)  }%
{\Gamma\left(  \tau-m\right)  }\left\{
\begin{array}
[c]{c}%
\left(
%TCIMACRO{\dint \limits_{0}^{1}}%
%BeginExpansion
{\displaystyle\int\limits_{0}^{1}}
%EndExpansion
\left\vert \frac{1}{3}\beta\left(  m+1,\tau-m\right)  -\frac{1}{2}\beta
_{w}\left(  m+1,\tau-m\right)  \right\vert ^{p}dw\right)  ^{\frac{1}{p}}\\
\times\left[
\begin{array}
[c]{c}%
\left(
%TCIMACRO{\dint \limits_{0}^{1}}%
%BeginExpansion
{\displaystyle\int\limits_{0}^{1}}
%EndExpansion
\left\vert
\begin{array}
[c]{c}%
\psi^{\prime}\left(  \frac{1+w}{2}\gamma+\frac{1-w}{2}\delta\right)
\end{array}
\right\vert ^{q}dw\right)  ^{\frac{1}{q}}\\
+\left(
%TCIMACRO{\dint \limits_{0}^{1}}%
%BeginExpansion
{\displaystyle\int\limits_{0}^{1}}
%EndExpansion
\left\vert
\begin{array}
[c]{c}%
\psi^{\prime}\left(  \frac{1-w}{2}\gamma+\frac{1+w}{2}\delta\right)
\end{array}
\right\vert ^{q}dw\right)  ^{\frac{1}{q}}%
\end{array}
\right]
\end{array}
\right\}  \\
&  \leq\frac{\delta-\gamma}{2.m!}\frac{\Gamma\left(  \tau+1\right)  }%
{\Gamma\left(  \tau-m\right)  }\left\{
\begin{array}
[c]{c}%
\left(
%TCIMACRO{\dint \limits_{0}^{1}}%
%BeginExpansion
{\displaystyle\int\limits_{0}^{1}}
%EndExpansion
\left\vert \frac{1}{3}\beta\left(  m+1,\tau-m\right)  -\frac{1}{2}\beta
_{w}\left(  m+1,\tau-m\right)  \right\vert ^{p}dt\right)  ^{\frac{1}{p}}\\
\times\left[
\begin{array}
[c]{c}%
\left(  \left\vert
\begin{array}
[c]{c}%
\psi^{\prime}\left(  \gamma\right)
\end{array}
\right\vert ^{q}%
%TCIMACRO{\dint \limits_{0}^{1}}%
%BeginExpansion
{\displaystyle\int\limits_{0}^{1}}
%EndExpansion
\frac{1+w}{2}dw+\left\vert
\begin{array}
[c]{c}%
\psi^{\prime}\left(  \delta\right)
\end{array}
\right\vert ^{q}%
%TCIMACRO{\dint \limits_{0}^{1}}%
%BeginExpansion
{\displaystyle\int\limits_{0}^{1}}
%EndExpansion
\frac{1-w}{2}dw\right)  ^{\frac{1}{q}}\\
+\left(  \left\vert
\begin{array}
[c]{c}%
\psi^{\prime}\left(  \gamma\right)
\end{array}
\right\vert ^{q}%
%TCIMACRO{\dint \limits_{0}^{1}}%
%BeginExpansion
{\displaystyle\int\limits_{0}^{1}}
%EndExpansion
\frac{1-w}{2}dw+\left\vert
\begin{array}
[c]{c}%
\psi^{\prime}\left(  \delta\right)
\end{array}
\right\vert ^{q}%
%TCIMACRO{\dint \limits_{0}^{1}}%
%BeginExpansion
{\displaystyle\int\limits_{0}^{1}}
%EndExpansion
\frac{1+w}{2}dw\right)  ^{\frac{1}{q}}%
\end{array}
\right]
\end{array}
\right\}  \\
&  \leq\frac{\delta-\gamma}{m!}\frac{\Gamma\left(  \tau+1\right)  }%
{\Gamma\left(  \tau-m\right)  }\left(  \frac{1}{2}\right)  ^{2q+1}\left(
Z_{2}\left(  \tau,m\right)  \right)  ^{\frac{1}{p}}\left[
\begin{array}
[c]{c}%
\left(  3\left\vert
\begin{array}
[c]{c}%
\psi^{\prime}\left(  \gamma\right)
\end{array}
\right\vert ^{q}+\left\vert
\begin{array}
[c]{c}%
\psi^{\prime}\left(  \delta\right)
\end{array}
\right\vert ^{q}\right)  ^{\frac{1}{q}}\\
+\left(  \left\vert
\begin{array}
[c]{c}%
\psi^{\prime}\left(  \gamma\right)
\end{array}
\right\vert ^{q}+3\left\vert
\begin{array}
[c]{c}%
\psi^{\prime}\left(  \delta\right)
\end{array}
\right\vert ^{q}\right)  ^{\frac{1}{q}}%
\end{array}
\right]  .
\end{align*}
This completes the proof.
\end{proof}

\begin{remark}
\label{3.6a} If we take $\tau=m+1,$ after that if we take $\tau=1$ in Theorem
\ref{3.6}, we obtain the following inequality
\begin{align*}
&  \left\vert \frac{1}{6}\left[  \psi\left(  \gamma\right)  +4\psi\left(
\frac{\gamma+\delta}{2}\right)  +\psi\left(  \delta\right)  \right]  -\frac
{1}{\delta-\gamma}%
%TCIMACRO{\dint \limits_{\gamma}^{\delta}}%
%BeginExpansion
{\displaystyle\int\limits_{\gamma}^{\delta}}
%EndExpansion
\psi\left(  \varepsilon\right)
\begin{array}
[c]{c}%
d\varepsilon
\end{array}
\right\vert \\
&  \leq\frac{\delta-\gamma}{12}\left(  \frac{2^{p+1}+1}{3\left(  p+1\right)
}\right)  ^{\frac{1}{p}}\left(  \frac{1}{4}\right)  ^{\frac{1}{q}}\\
&  \times\left[
\begin{array}
[c]{c}%
\left(  3\left\vert
\begin{array}
[c]{c}%
\psi^{\prime}\left(  \gamma\right)
\end{array}
\right\vert ^{q}+\left\vert
\begin{array}
[c]{c}%
\psi^{\prime}\left(  \delta\right)
\end{array}
\right\vert ^{q}\right)  ^{\frac{1}{q}}\\
+\left(  \left\vert
\begin{array}
[c]{c}%
\psi^{\prime}\left(  \gamma\right)
\end{array}
\right\vert ^{q}+3\left\vert
\begin{array}
[c]{c}%
\psi^{\prime}\left(  \delta\right)
\end{array}
\right\vert ^{q}\right)  ^{\frac{1}{q}}%
\end{array}
\right]
\end{align*}
with%
\[%
%TCIMACRO{\dint \limits_{0}^{1}}%
%BeginExpansion
{\displaystyle\int\limits_{0}^{1}}
%EndExpansion
\left\vert \frac{1}{2}w-\frac{1}{3}\right\vert ^{p}dw=\frac{2^{p+2}+2}{\left(
p+1\right)  6^{p+1}}.
\]

\end{remark}

See also (\cite{SSO10}, Theorem 4).

\begin{theorem}
\label{3.7}Let $\psi:I\subset\left(  0,\infty\right)  \rightarrow%
%TCIMACRO{\U{211d} }%
%BeginExpansion
\mathbb{R}
%EndExpansion
,$ be a differentiable function on $I^{\circ},\gamma,\delta\in I^{\circ}$ and
$\gamma<\delta.$ If $\psi^{\prime}$ $\in L\left[  \gamma,\delta\right]  $ and
$\left\vert \psi^{\prime}\right\vert ^{q}$ is a convex function on $\left[
\gamma,\delta\right]  $ for $q\geq1,$ then the following inequality holds:
\begin{align}
&  \left\vert \frac{1}{6}\left[  \psi\left(  \gamma\right)  +4\psi\left(
\frac{\gamma+\delta}{2}\right)  +\psi\left(  \delta\right)  \right]
-\frac{2^{\tau-1}}{\left(  \delta-\gamma\right)  ^{\tau}}\frac{\Gamma\left(
\tau+1\right)  }{\Gamma\left(  \tau-m\right)  }\left[
\begin{array}
[c]{c}%
I_{\tau}^{\gamma}\psi\left(  \frac{\gamma+\delta}{2}\right)
\end{array}
+%
\begin{array}
[c]{c}%
^{\delta}I_{\tau}\psi\left(  \frac{\gamma+\delta}{2}\right)
\end{array}
\right]  \right\vert \label{3-7}\\
&  \leq\frac{\delta-\gamma}{2.m!}\frac{\Gamma\left(  \tau+1\right)  }%
{\Gamma\left(  \tau-m\right)  }\left(  Z_{1}\left(  \tau,m\right)  \right)
^{1-\frac{1}{q}}\left[
\begin{array}
[c]{c}%
\left(  Z_{3}\left(  \tau,m\right)  \left\vert
\begin{array}
[c]{c}%
\psi^{\prime}\left(  \gamma\right)
\end{array}
\right\vert ^{q}+Z_{4}\left(  \tau,m\right)  \left\vert
\begin{array}
[c]{c}%
\psi^{\prime}\left(  \delta\right)
\end{array}
\right\vert ^{q}\right)  ^{\frac{1}{q}}\\
+\left(  Z_{4}\left(  \tau,m\right)  \left\vert
\begin{array}
[c]{c}%
\psi^{\prime}\left(  \gamma\right)
\end{array}
\right\vert ^{q}+Z_{3}\left(  \tau,m\right)  \left\vert
\begin{array}
[c]{c}%
\psi^{\prime}\left(  \delta\right)
\end{array}
\right\vert ^{q}\right)  ^{\frac{1}{q}}%
\end{array}
\right] \nonumber
\end{align}
where%
\[
Z_{3}\left(  \tau,m\right)  =%
%TCIMACRO{\dint \limits_{0}^{1}}%
%BeginExpansion
{\displaystyle\int\limits_{0}^{1}}
%EndExpansion
\left\vert \frac{1}{3}\beta\left(  m+1,\tau-m\right)  -\frac{1}{2}\beta
_{w}\left(  m+1,\tau-m\right)  \right\vert \frac{1+w}{2}dw
\]%
\[
Z_{4}\left(  \tau,m\right)  =%
%TCIMACRO{\dint \limits_{0}^{1}}%
%BeginExpansion
{\displaystyle\int\limits_{0}^{1}}
%EndExpansion
\left\vert \frac{1}{3}\beta\left(  m+1,\tau-m\right)  -\frac{1}{2}\beta
_{w}\left(  m+1,\tau-m\right)  \right\vert \frac{1-w}{2}dw
\]

\end{theorem}

and $Z_{1}\left(  \tau,m\right)  $ is defined as in the Theorem \ref{3.3} with
$m=0,1,2,...$ and $\tau\in\left(  m,m+1\right]  .$

\begin{proof}
From Lemma \ref{3.1} and using the power mean inequality, we have that the
following inequality holds:
\begin{align*}
&  \left\vert \frac{1}{6}\left[  \psi\left(  \gamma\right)  +4\psi\left(
\frac{\gamma+\delta}{2}\right)  +\psi\left(  \delta\right)  \right]
-\frac{2^{\tau-1}}{\left(  \delta-\gamma\right)  ^{\tau}}\frac{\Gamma\left(
\tau+1\right)  }{\Gamma\left(  \tau-m\right)  }\left[
\begin{array}
[c]{c}%
I_{\tau}^{\gamma}\psi\left(  \frac{\gamma+\delta}{2}\right)
\end{array}
+%
\begin{array}
[c]{c}%
^{\delta}I_{\tau}\psi\left(  \frac{\gamma+\delta}{2}\right)
\end{array}
\right]  \right\vert \\
&  \leq\frac{\delta-\gamma}{2.m!}\frac{\Gamma\left(  \tau+1\right)  }%
{\Gamma\left(  \tau-m\right)  }\left\{
\begin{array}
[c]{c}%
%TCIMACRO{\dint \limits_{0}^{1}}%
%BeginExpansion
{\displaystyle\int\limits_{0}^{1}}
%EndExpansion
\left\vert \frac{1}{3}\beta\left(  m+1,\tau-m\right)  -\frac{1}{2}\beta
_{w}\left(  m+1,\tau-m\right)  \right\vert \\
\times\left\vert
\begin{array}
[c]{c}%
\psi^{\prime}\left(  \frac{1+w}{2}\gamma+\frac{1-w}{2}\delta\right)
\end{array}
\right\vert dw\\
+%
%TCIMACRO{\dint \limits_{0}^{1}}%
%BeginExpansion
{\displaystyle\int\limits_{0}^{1}}
%EndExpansion
\left\vert \frac{1}{3}\beta\left(  m+1,\tau-m\right)  -\frac{1}{2}\beta
_{w}\left(  m+1,\tau-m\right)  \right\vert \\
\times\left\vert
\begin{array}
[c]{c}%
\psi^{\prime}\left(  \frac{1-w}{2}\gamma+\frac{1+w}{2}\delta\right)
\end{array}
\right\vert dw
\end{array}
\right\}  \\
&  \leq\frac{\delta-\gamma}{2.m!}\frac{\Gamma\left(  \tau+1\right)  }%
{\Gamma\left(  \tau-m\right)  }\left\{
\begin{array}
[c]{c}%
\left(
%TCIMACRO{\dint \limits_{0}^{1}}%
%BeginExpansion
{\displaystyle\int\limits_{0}^{1}}
%EndExpansion
\left\vert \frac{1}{3}\beta\left(  m+1,\tau-m\right)  -\frac{1}{2}\beta
_{w}\left(  m+1,\tau-m\right)  \right\vert dw\right)  ^{1-\frac{1}{q}}\\
\times\left[
\begin{array}
[c]{c}%
\left(
\begin{array}
[c]{c}%
%TCIMACRO{\dint \limits_{0}^{1}}%
%BeginExpansion
{\displaystyle\int\limits_{0}^{1}}
%EndExpansion
\left\vert \frac{1}{3}\beta\left(  m+1,\tau-m\right)  -\frac{1}{2}\beta
_{w}\left(  m+1,\tau-m\right)  \right\vert \\
\times\left\vert
\begin{array}
[c]{c}%
\psi^{\prime}\left(  \frac{1+w}{2}\gamma+\frac{1-w}{2}\delta\right)
\end{array}
\right\vert ^{q}dw
\end{array}
\right)  ^{\frac{1}{q}}\\
+\left(
\begin{array}
[c]{c}%
%TCIMACRO{\dint \limits_{0}^{1}}%
%BeginExpansion
{\displaystyle\int\limits_{0}^{1}}
%EndExpansion
\left\vert \frac{1}{3}\beta\left(  m+1,\tau-m\right)  -\frac{1}{2}\beta
_{w}\left(  m+1,\tau-m\right)  \right\vert \\
\times\left\vert
\begin{array}
[c]{c}%
\psi^{\prime}\left(  \frac{1-w}{2}\gamma+\frac{1+w}{2}\delta\right)
\end{array}
\right\vert ^{q}dw
\end{array}
\right)  ^{\frac{1}{q}}%
\end{array}
\right]
\end{array}
\right\}  .
\end{align*}

By the convexity of $\left\vert \psi^{\prime}\right\vert ^{q}$ \ %

\begin{align*}
&
%TCIMACRO{\dint \limits_{0}^{1}}%
%BeginExpansion
{\displaystyle\int\limits_{0}^{1}}
%EndExpansion
\left\vert \frac{1}{3}\beta\left(  m+1,\tau-m\right)  -\frac{1}{2}\beta
_{w}\left(  m+1,\tau-m\right)  \right\vert \left\vert
\begin{array}
[c]{c}%
\psi^{\prime}\left(  \frac{1+w}{2}\gamma+\frac{1-w}{2}\delta\right)
\end{array}
\right\vert ^{q}dw\\
&  \leq\left\vert
\begin{array}
[c]{c}%
\psi^{\prime}\left(  \gamma\right)
\end{array}
\right\vert ^{q}%
%TCIMACRO{\dint \limits_{0}^{1}}%
%BeginExpansion
{\displaystyle\int\limits_{0}^{1}}
%EndExpansion
\left\vert \frac{1}{3}\beta\left(  m+1,\tau-m\right)  -\frac{1}{2}\beta
_{w}\left(  m+1,\tau-m\right)  \right\vert \frac{1+w}{2}dw\\
&  +\left\vert
\begin{array}
[c]{c}%
\psi^{\prime}\left(  \delta\right)
\end{array}
\right\vert ^{q}%
%TCIMACRO{\dint \limits_{0}^{1}}%
%BeginExpansion
{\displaystyle\int\limits_{0}^{1}}
%EndExpansion
\left\vert \frac{1}{3}\beta\left(  m+1,\tau-m\right)  -\frac{1}{2}\beta
_{w}\left(  m+1,\tau-m\right)  \right\vert \frac{1-w}{2}dw
\end{align*}

and%
\begin{align*}
&
%TCIMACRO{\dint \limits_{0}^{1}}%
%BeginExpansion
{\displaystyle\int\limits_{0}^{1}}
%EndExpansion
\left\vert \frac{1}{3}\beta\left(  m+1,\tau-m\right)  -\frac{1}{2}\beta
_{w}\left(  m+1,\tau-m\right)  \right\vert \left\vert
\begin{array}
[c]{c}%
\psi^{\prime}\left(  \frac{1-w}{2}\gamma+\frac{1+w}{2}\delta\right)
\end{array}
\right\vert ^{q}dt\\
&  \leq\left\vert
\begin{array}
[c]{c}%
\psi^{\prime}\left(  \gamma\right)
\end{array}
\right\vert ^{q}%
%TCIMACRO{\dint \limits_{0}^{1}}%
%BeginExpansion
{\displaystyle\int\limits_{0}^{1}}
%EndExpansion
\left\vert \frac{1}{3}\beta\left(  m+1,\tau-m\right)  -\frac{1}{2}\beta
_{w}\left(  m+1,\tau-m\right)  \right\vert \frac{1-w}{2}dw\\
&  +\left\vert
\begin{array}
[c]{c}%
\psi^{\prime}\left(  \delta\right)
\end{array}
\right\vert ^{q}%
%TCIMACRO{\dint \limits_{0}^{1}}%
%BeginExpansion
{\displaystyle\int\limits_{0}^{1}}
%EndExpansion
\left\vert \frac{1}{3}\beta\left(  m+1,\tau-m\right)  -\frac{1}{2}\beta
_{w}\left(  m+1,\tau-m\right)  \right\vert \frac{1+w}{2}dw
\end{align*}

Using the last two inequalities we obtain%

\begin{align*}
&  \left\vert \frac{1}{6}\left[  \psi\left(  \gamma\right)  +4\psi\left(
\frac{\gamma+\delta}{2}\right)  +\psi\left(  \delta\right)  \right]
-\frac{2^{\tau-1}}{\left(  \delta-\gamma\right)  ^{\tau}}\frac{\Gamma\left(
\tau+1\right)  }{\Gamma\left(  \tau-m\right)  }\left[
\begin{array}
[c]{c}%
I_{\tau}^{\gamma}\psi\left(  \frac{\gamma+\delta}{2}\right)
\end{array}
+%
\begin{array}
[c]{c}%
^{\delta}I_{\tau}\psi\left(  \frac{\gamma+\delta}{2}\right)
\end{array}
\right]  \right\vert \\
&  \leq\frac{\delta-\gamma}{2.m!}\frac{\Gamma\left(  \tau+1\right)  }%
{\Gamma\left(  \tau-m\right)  }\left(
%TCIMACRO{\dint \limits_{0}^{1}}%
%BeginExpansion
{\displaystyle\int\limits_{0}^{1}}
%EndExpansion
\left\vert \frac{1}{3}\beta\left(  m+1,\tau-m\right)  -\frac{1}{2}\beta
_{w}\left(  m+1,\tau-m\right)  \right\vert dw\right)  ^{1-\frac{1}{q}}\\
&  \times\left[
\begin{array}
[c]{c}%
\left(
\begin{array}
[c]{c}%
\left\vert
\begin{array}
[c]{c}%
\psi^{\prime}\left(  \gamma\right)
\end{array}
\right\vert ^{q}%
%TCIMACRO{\dint \limits_{0}^{1}}%
%BeginExpansion
{\displaystyle\int\limits_{0}^{1}}
%EndExpansion
\left\vert \frac{1}{3}\beta\left(  m+1,\tau-m\right)  -\frac{1}{2}\beta
_{w}\left(  m+1,\tau-m\right)  \right\vert \frac{1+w}{2}dw\\
+\left\vert
\begin{array}
[c]{c}%
\psi^{\prime}\left(  \delta\right)
\end{array}
\right\vert ^{q}%
%TCIMACRO{\dint \limits_{0}^{1}}%
%BeginExpansion
{\displaystyle\int\limits_{0}^{1}}
%EndExpansion
\left\vert \frac{1}{3}\beta\left(  m+1,\tau-m\right)  -\frac{1}{2}\beta
_{w}\left(  m+1,\tau-m\right)  \right\vert \frac{1-w}{2}dw
\end{array}
\right)  ^{\frac{1}{q}}\\
+\left(
\begin{array}
[c]{c}%
\left\vert
\begin{array}
[c]{c}%
\psi^{\prime}\left(  \gamma\right)
\end{array}
\right\vert ^{q}%
%TCIMACRO{\dint \limits_{0}^{1}}%
%BeginExpansion
{\displaystyle\int\limits_{0}^{1}}
%EndExpansion
\left\vert \frac{1}{3}\beta\left(  m+1,\tau-m\right)  -\frac{1}{2}\beta
_{w}\left(  m+1,\tau-m\right)  \right\vert \frac{1-w}{2}dw\\
+\left\vert
\begin{array}
[c]{c}%
\psi^{\prime}\left(  \delta\right)
\end{array}
\right\vert ^{q}%
%TCIMACRO{\dint \limits_{0}^{1}}%
%BeginExpansion
{\displaystyle\int\limits_{0}^{1}}
%EndExpansion
\left\vert \frac{1}{3}\beta\left(  m+1,\tau-m\right)  -\frac{1}{2}\beta
_{w}\left(  m+1,\tau-m\right)  \right\vert \frac{1+w}{2}dw
\end{array}
\right)  ^{\frac{1}{q}}%
\end{array}
\right]  .
\end{align*}

\end{proof}

\begin{remark}
\label{3.8} If we take $\tau=m+1,$ after that if we take $\tau=1$ in Theorem
\ref{3.7}, we obtain the following inequality
\begin{align}
&  \left\vert \frac{1}{6}\left[  \psi\left(  \gamma\right)  +4\psi\left(
\frac{\gamma+\delta}{2}\right)  +\psi\left(  \delta\right)  \right]  -\frac
{1}{\delta-\gamma}%
%TCIMACRO{\dint \limits_{\gamma}^{\delta}}%
%BeginExpansion
{\displaystyle\int\limits_{\gamma}^{\delta}}
%EndExpansion
\psi\left(  \varepsilon\right)  d\varepsilon\right\vert \label{3-8}\\
&  \leq\frac{\delta-\gamma}{2}\left(  \frac{5}{36}\right)  ^{1-\frac{1}{q}%
}\left(  \frac{1}{648}\right)  ^{\frac{1}{q}}\left[
\begin{array}
[c]{c}%
\left(  61\left\vert
\begin{array}
[c]{c}%
\psi^{\prime}\left(  \gamma\right)
\end{array}
\right\vert ^{q}+29\left\vert
\begin{array}
[c]{c}%
\psi^{\prime}\left(  \delta\right)
\end{array}
\right\vert ^{q}\right)  ^{\frac{1}{q}}\\
+\left(  29\left\vert
\begin{array}
[c]{c}%
\psi^{\prime}\left(  \gamma\right)
\end{array}
\right\vert ^{q}+61\left\vert
\begin{array}
[c]{c}%
\psi^{\prime}\left(  \delta\right)
\end{array}
\right\vert ^{q}\right)  ^{\frac{1}{q}}%
\end{array}
\right]  .\nonumber
\end{align}

\end{remark}

See also (\cite{SB19}, Remark 2.7.).

\begin{theorem}
\label{3.9}Let $\psi:I\subset\left(  0,\infty\right)  \rightarrow%
%TCIMACRO{\U{211d} }%
%BeginExpansion
\mathbb{R}
%EndExpansion
,$ be a differentiable function on $I^{\circ},\gamma,\delta\in I^{\circ}$ and
$\gamma<\delta.$ If $f^{\prime}$ $\in L\left(  \left[  \gamma,\delta\right]
\right)  $ and $\left\vert \psi^{\prime}\right\vert ^{q}$ is a convex function
on $\left[  \gamma,\delta\right]  $ for $q>1$ and $\frac{1}{p}+\frac{1}{q}=1,$
then the following inequality holds:
\begin{align}
&  \left\vert \frac{1}{6}\left[  \psi\left(  \gamma\right)  +4\psi\left(
\frac{\gamma+\delta}{2}\right)  +\psi\left(  \delta\right)  \right]
-\frac{2^{\tau-1}}{\left(  \delta-\gamma\right)  ^{\tau}}\frac{\Gamma\left(
\tau+1\right)  }{\Gamma\left(  \tau-m\right)  }\left[
\begin{array}
[c]{c}%
I_{\tau}^{\gamma}\psi\left(  \frac{\gamma+\delta}{2}\right)
\end{array}
+%
\begin{array}
[c]{c}%
^{\delta}I_{\tau}\psi\left(  \frac{\gamma+\delta}{2}\right)
\end{array}
\right]  \right\vert \label{3-9}\\
&  \leq\frac{\delta-\gamma}{2.m!}\frac{\Gamma\left(  \tau+1\right)  }%
{\Gamma\left(  \tau-m\right)  }\left(  \frac{1}{2}\right)  ^{\frac{1}{q}%
}\left(  Z_{2}\left(  \tau,m\right)  \right)  ^{\frac{1}{p}}\left[
\begin{array}
[c]{c}%
\left(  \left\vert \psi^{\prime}\left(  \gamma\right)  \right\vert
^{q}+\left\vert \psi^{\prime}\left(  \frac{\gamma+\delta}{2}\right)
\right\vert ^{q}\right)  ^{\frac{1}{q}}\\
+\left(  \left\vert \psi^{\prime}\left(  \frac{\gamma+\delta}{2}\right)
\right\vert ^{q}+\left\vert \psi^{\prime}\left(  \delta\right)  \right\vert
^{q}\right)  ^{\frac{1}{q}}%
\end{array}
\right]  ,\nonumber
\end{align}
where $Z_{2}\left(  \tau,m\right)  $ is defined as in Theorem \ref{3.6} with
$m=0,1,2,...$ and $\tau\in\left(  m,m+1\right]  .$
\end{theorem}

\begin{proof}
From Lemma \ref{3.1} and using the H\"{o}lder's inequality, we have
\begin{align*}
&  \left\vert \frac{1}{6}\left[  \psi\left(  \gamma\right)  +4\psi\left(
\frac{\gamma+\delta}{2}\right)  +\psi\left(  \delta\right)  \right]
-\frac{2^{\tau-1}}{\left(  \delta-\gamma\right)  ^{\tau}}\frac{\Gamma\left(
\tau+1\right)  }{\Gamma\left(  \tau-m\right)  }\left[
\begin{array}
[c]{c}%
I_{\tau}^{\gamma}\psi\left(  \frac{\gamma+\delta}{2}\right)
\end{array}
+%
\begin{array}
[c]{c}%
^{\delta}I_{\tau}\psi\left(  \frac{\gamma+\delta}{2}\right)
\end{array}
\right]  \right\vert \\
&  \leq\frac{\delta-\gamma}{2.m!}\frac{\Gamma\left(  \tau+1\right)  }%
{\Gamma\left(  \tau-m\right)  }\left\{
\begin{array}
[c]{c}%
\left(
%TCIMACRO{\dint \limits_{0}^{1}}%
%BeginExpansion
{\displaystyle\int\limits_{0}^{1}}
%EndExpansion
\left\vert \frac{1}{3}\beta\left(  m+1,\tau-m\right)  -\frac{1}{2}\beta
_{w}\left(  m+1,\tau-m\right)  \right\vert ^{p}dw\right)  ^{\frac{1}{p}}\\
\times\left(
%TCIMACRO{\dint \limits_{0}^{1}}%
%BeginExpansion
{\displaystyle\int\limits_{0}^{1}}
%EndExpansion
\left\vert
\begin{array}
[c]{c}%
\psi^{\prime}\left(  \frac{1+w}{2}\gamma+\frac{1-w}{2}\delta\right)
\end{array}
\right\vert ^{q}dw\right)  ^{\frac{1}{q}}\\
+\left(
%TCIMACRO{\dint \limits_{0}^{1}}%
%BeginExpansion
{\displaystyle\int\limits_{0}^{1}}
%EndExpansion
\left\vert \frac{1}{3}\beta\left(  m+1,\tau-m\right)  -\frac{1}{2}\beta
_{w}\left(  m+1,\tau-m\right)  \right\vert ^{p}dw\right)  ^{\frac{1}{p}}\\
\times\left(
%TCIMACRO{\dint \limits_{0}^{1}}%
%BeginExpansion
{\displaystyle\int\limits_{0}^{1}}
%EndExpansion
\left\vert
\begin{array}
[c]{c}%
\psi^{\prime}\left(  \frac{1-w}{2}\gamma+\frac{1+w}{2}\delta\right)
\end{array}
\right\vert ^{q}dw\right)  ^{\frac{1}{q}}%
\end{array}
\right\}  \\
&  \leq\frac{\delta-\gamma}{2.m!}\frac{\Gamma\left(  \tau+1\right)  }%
{\Gamma\left(  \tau-m\right)  }\left\{
\begin{array}
[c]{c}%
\left(
%TCIMACRO{\dint \limits_{0}^{1}}%
%BeginExpansion
{\displaystyle\int\limits_{0}^{1}}
%EndExpansion
\left\vert \frac{1}{3}\beta\left(  m+1,\tau-m\right)  -\frac{1}{2}\beta
_{w}\left(  m+1,\tau-m\right)  \right\vert ^{p}dw\right)  ^{\frac{1}{p}}\\
\times\left[
\begin{array}
[c]{c}%
\left(
%TCIMACRO{\dint \limits_{0}^{1}}%
%BeginExpansion
{\displaystyle\int\limits_{0}^{1}}
%EndExpansion
\left\vert
\begin{array}
[c]{c}%
\psi^{\prime}\left(  \frac{1+w}{2}\gamma+\frac{1-w}{2}\delta\right)
\end{array}
\right\vert ^{q}dw\right)  ^{\frac{1}{q}}\\
+\left(
%TCIMACRO{\dint \limits_{0}^{1}}%
%BeginExpansion
{\displaystyle\int\limits_{0}^{1}}
%EndExpansion
\left\vert
\begin{array}
[c]{c}%
\psi^{\prime}\left(  \frac{1-w}{2}\gamma+\frac{1+w}{2}\delta\right)
\end{array}
\right\vert ^{q}dw\right)  ^{\frac{1}{q}}%
\end{array}
\right]
\end{array}
\right\}  .
\end{align*}

Since $\left\vert \psi^{\prime}\right\vert ^{q}$ is covex by \ \ref{1-3} we have%

\[%
%TCIMACRO{\dint \limits_{0}^{1}}%
%BeginExpansion
{\displaystyle\int\limits_{0}^{1}}
%EndExpansion
\left\vert
\begin{array}
[c]{c}%
\psi^{\prime}\left(  \frac{1+w}{2}\gamma+\frac{1-w}{2}\delta\right)
\end{array}
\right\vert ^{q}dw\leq\frac{\left\vert \psi^{\prime}\left(  \gamma\right)
\right\vert ^{q}+\left\vert \psi^{\prime}\left(  \frac{\gamma+\delta}%
{2}\right)  \right\vert ^{q}}{2},
\]

and%
\[%
%TCIMACRO{\dint \limits_{0}^{1}}%
%BeginExpansion
{\displaystyle\int\limits_{0}^{1}}
%EndExpansion
\left\vert
\begin{array}
[c]{c}%
\psi^{\prime}\left(  \frac{1-w}{2}\gamma+\frac{1+w}{2}\delta\right)
\end{array}
\right\vert ^{q}dw\leq\frac{\left\vert \psi^{\prime}\left(  \frac
{\gamma+\delta}{2}\right)  \right\vert ^{q}+\left\vert \psi^{\prime}\left(
\delta\right)  \right\vert ^{q}}{2}.
\]

So, we obtain%

\begin{align*}
&  \left\vert \frac{1}{6}\left[  \psi\left(  \gamma\right)  +4\psi\left(
\frac{\gamma+\delta}{2}\right)  +\psi\left(  \delta\right)  \right]
-\frac{2^{\tau-1}}{\left(  \delta-\gamma\right)  ^{\tau}}\frac{\Gamma\left(
\tau+1\right)  }{\Gamma\left(  \tau-m\right)  }\left[
\begin{array}
[c]{c}%
I_{\tau}^{\gamma}\psi\left(  \frac{\gamma+\delta}{2}\right)
\end{array}
+%
\begin{array}
[c]{c}%
^{\delta}I_{\tau}\psi\left(  \frac{\gamma+\delta}{2}\right)
\end{array}
\right]  \right\vert \\
&  \leq\frac{\delta-\gamma}{2.m!}\frac{\Gamma\left(  \tau+1\right)  }%
{\Gamma\left(  \tau-m\right)  }\left(  \frac{1}{2}\right)  ^{\frac{1}{q}%
}\left\{
\begin{array}
[c]{c}%
\left(
%TCIMACRO{\dint \limits_{0}^{1}}%
%BeginExpansion
{\displaystyle\int\limits_{0}^{1}}
%EndExpansion
\left\vert \frac{1}{3}\beta\left(  m+1,\tau-m\right)  -\frac{1}{2}\beta
_{w}\left(  m+1,\tau-m\right)  \right\vert ^{p}dx\right)  ^{\frac{1}{p}}\\
\times\left[
\begin{array}
[c]{c}%
\left(  \left\vert \psi^{\prime}\left(  \gamma\right)  \right\vert
^{q}+\left\vert \psi^{\prime}\left(  \frac{\gamma+\delta}{2}\right)
\right\vert ^{q}\right)  ^{\frac{1}{q}}\\
+\left(  \left\vert \psi^{\prime}\left(  \frac{\gamma+\delta}{2}\right)
\right\vert ^{q}+\left\vert \psi^{\prime}\left(  \delta\right)  \right\vert
^{q}\right)  ^{\frac{1}{q}}%
\end{array}
\right]
\end{array}
\right\}  .
\end{align*}

\end{proof}

\begin{remark}
\label{3.10}If we take $\tau=m+1,$ after that if we take $\tau=1$ in Theorem
\ref{3.6}, we obtain the following inequality
\begin{align*}
&  \left\vert \frac{1}{6}\left[  \psi\left(  \gamma\right)  +4\psi\left(
\frac{\gamma+\delta}{2}\right)  +\psi\left(  \delta\right)  \right]  -\frac
{1}{\delta-\gamma}%
%TCIMACRO{\dint \limits_{\gamma}^{\delta}}%
%BeginExpansion
{\displaystyle\int\limits_{\gamma}^{\delta}}
%EndExpansion
\psi\left(  \varepsilon\right)
\begin{array}
[c]{c}%
d\varepsilon
\end{array}
\right\vert \\
&  \leq\frac{\delta-\gamma}{12}\left(  \frac{2^{p+1}+1}{3\left(  p+1\right)
}\right)  ^{\frac{1}{p}}\left(  \frac{1}{4}\right)  ^{\frac{1}{q}}\\
&  \times\left[
\begin{array}
[c]{c}%
\left(  \left\vert
\begin{array}
[c]{c}%
\psi^{\prime}\left(  \gamma\right)
\end{array}
\right\vert ^{q}+3\left\vert
\begin{array}
[c]{c}%
\psi^{\prime}\left(  \frac{\gamma+\delta}{2}\right)
\end{array}
\right\vert ^{q}\right)  ^{\frac{1}{q}}\\
+\left(  3\left\vert
\begin{array}
[c]{c}%
\psi^{\prime}\left(  \frac{\gamma+\delta}{2}\right)
\end{array}
\right\vert ^{q}+\left\vert
\begin{array}
[c]{c}%
\psi^{\prime}\left(  \delta\right)
\end{array}
\right\vert ^{q}\right)  ^{\frac{1}{q}}%
\end{array}
\right]
\end{align*}
with%
\[%
%TCIMACRO{\dint \limits_{0}^{1}}%
%BeginExpansion
{\displaystyle\int\limits_{0}^{1}}
%EndExpansion
\left\vert \frac{1}{2}w-\frac{1}{3}\right\vert ^{p}dw=\frac{2^{p+2}+2}{\left(
p+1\right)  6^{p+1}}.
\]

\end{remark}

See also (\cite{SSO10}, Theorem 9 (for $s=1$)).

\section{Estimation Results}

In the function $\psi^{\prime}$ is bounded from below and above, then we have
the following result.

\begin{theorem}
\label{4.1}Let $\psi:\left[  \gamma,\delta\right]  \rightarrow%
%TCIMACRO{\U{211d} }%
%BeginExpansion
\mathbb{R}
%EndExpansion
$ be differentiable and continuous mappings on $\left(  \gamma,\delta\right)
$ and let $\psi^{\prime}\in L\left[  \gamma,\delta\right]  .$ Assume that
there exist constants $k<K$ such that $-\infty<k\leq%
\begin{array}
[c]{c}%
\psi^{\prime}%
\end{array}
\leq K<+\infty.$ Then,%
\begin{align}
&  \left\vert
\begin{array}
[c]{c}%
\frac{1}{6}\left[  \psi\left(  \gamma\right)  +4\psi\left(  \frac
{\gamma+\delta}{2}\right)  +\psi\left(  \delta\right)  \right]  -\frac
{2^{\tau-1}}{\left(  \delta-\gamma\right)  ^{\tau}}\frac{\Gamma\left(
\tau+1\right)  }{\Gamma\left(  \tau-m\right)  }\left[
\begin{array}
[c]{c}%
I_{\tau}^{\gamma}\psi\left(  \frac{\gamma+\delta}{2}\right)
\end{array}
+%
\begin{array}
[c]{c}%
^{\delta}I_{\tau}\psi\left(  \frac{\gamma+\delta}{2}\right)
\end{array}
\right] \\
-\frac{K+k}{2}\left(  \delta-\gamma\right)
%TCIMACRO{\dint \limits_{0}^{1}}%
%BeginExpansion
{\displaystyle\int\limits_{0}^{1}}
%EndExpansion
h\left(  w\right)  dw
\end{array}
\right\vert \label{4-1}\\
&  \leq\frac{K-k}{2.m!}\frac{\Gamma\left(  \tau+1\right)  }{\Gamma\left(
\tau-m\right)  }\left(  \delta-\gamma\right)  Z_{1}\left(  \tau,m\right)
,\nonumber
\end{align}

\end{theorem}

where
\[
h\left(  w\right)  =\frac{1}{3}%
%TCIMACRO{\dint \limits_{0}^{1}}%
%BeginExpansion
{\displaystyle\int\limits_{0}^{1}}
%EndExpansion
\beta\left(  m+1,\tau-m\right)  dw-\frac{1}{2}%
%TCIMACRO{\dint \limits_{0}^{1}}%
%BeginExpansion
{\displaystyle\int\limits_{0}^{1}}
%EndExpansion
\beta_{w}\left(  m+1,\tau-m\right)  dw.
\]

\begin{proof}
From Lemma \ref{3.1}, we have that%
\begin{align*}
&  \frac{1}{6}\left[  \psi\left(  \gamma\right)  +4\psi\left(  \frac
{\gamma+\delta}{2}\right)  +\psi\left(  \delta\right)  \right]  -\frac
{2^{\tau-1}}{\left(  \delta-\gamma\right)  ^{\tau}}\frac{\Gamma\left(
\tau+1\right)  }{\Gamma\left(  \tau-m\right)  }\left[
\begin{array}
[c]{c}%
I_{\tau}^{\gamma}\psi\left(  \frac{\gamma+\delta}{2}\right)
\end{array}
+%
\begin{array}
[c]{c}%
^{\delta}I_{\tau}\psi\left(  \frac{\gamma+\delta}{2}\right)
\end{array}
\right] \\
&  =\frac{\delta-\gamma}{2.m!}\frac{\Gamma\left(  \tau+1\right)  }%
{\Gamma\left(  \tau-m\right)  }\left\{
\begin{array}
[c]{c}%
%TCIMACRO{\dint \limits_{0}^{1}}%
%BeginExpansion
{\displaystyle\int\limits_{0}^{1}}
%EndExpansion
h\left(  w\right)  \left[
\begin{array}
[c]{c}%
\psi^{\prime}\left(  \frac{1+w}{2}\gamma+\frac{1-w}{2}\delta\right)
\end{array}
-\frac{K+k}{2}+\frac{K+k}{2}\right]  dw\\
+%
%TCIMACRO{\dint \limits_{0}^{1}}%
%BeginExpansion
{\displaystyle\int\limits_{0}^{1}}
%EndExpansion
h\left(  w\right)  \left[
\begin{array}
[c]{c}%
\psi^{\prime}\left(  \frac{1-w}{2}\gamma+\frac{1+w}{2}\delta\right)
\end{array}
-\frac{K+k}{2}+\frac{K+k}{2}\right]  dw
\end{array}
\right\} \\
&  =\frac{\delta-\gamma}{2.m!}\frac{\Gamma\left(  \tau+1\right)  }%
{\Gamma\left(  \tau-m\right)  }\left\{
\begin{array}
[c]{c}%
%TCIMACRO{\dint \limits_{0}^{1}}%
%BeginExpansion
{\displaystyle\int\limits_{0}^{1}}
%EndExpansion
h\left(  w\right)  \left[
\begin{array}
[c]{c}%
\psi^{\prime}\left(  \frac{1+w}{2}\gamma+\frac{1-w}{2}\delta\right)
\end{array}
-\frac{K+k}{2}\right]  dw\\
+%
%TCIMACRO{\dint \limits_{0}^{1}}%
%BeginExpansion
{\displaystyle\int\limits_{0}^{1}}
%EndExpansion
h\left(  w\right)  \left[
\begin{array}
[c]{c}%
\psi^{\prime}\left(  \frac{1-w}{2}\gamma+\frac{1+w}{2}\delta\right)
\end{array}
-\frac{K+k}{2}\right]  dw
\end{array}
\right\} \\
&  +\frac{K+k}{2}\left(  \delta-\gamma\right)
%TCIMACRO{\dint \limits_{0}^{1}}%
%BeginExpansion
{\displaystyle\int\limits_{0}^{1}}
%EndExpansion
h\left(  w\right)  dw.
\end{align*}

So%
\begin{align*}
&  \left\vert
\begin{array}
[c]{c}%
\frac{1}{6}\left[  \psi\left(  \gamma\right)  +4\psi\left(  \frac
{\gamma+\delta}{2}\right)  +\psi\left(  \delta\right)  \right]  -\frac
{2^{\tau-1}}{\left(  \delta-\gamma\right)  ^{\tau}}\frac{\Gamma\left(
\tau+1\right)  }{\Gamma\left(  \tau-m\right)  }\left[
\begin{array}
[c]{c}%
I_{\tau}^{\gamma}\psi\left(  \frac{\gamma+\delta}{2}\right)
\end{array}
+%
\begin{array}
[c]{c}%
^{\delta}I_{\tau}\psi\left(  \frac{\gamma+\delta}{2}\right)
\end{array}
\right] \\
-\frac{K+k}{2}\left(  \delta-\gamma\right)
%TCIMACRO{\dint \limits_{0}^{1}}%
%BeginExpansion
{\displaystyle\int\limits_{0}^{1}}
%EndExpansion
h\left(  w\right)  dw
\end{array}
\right\vert \\
&  \leq\frac{\delta-\gamma}{2.m!}\frac{\Gamma\left(  \tau+1\right)  }%
{\Gamma\left(  \tau-m\right)  }\left\{
\begin{array}
[c]{c}%
%TCIMACRO{\dint \limits_{0}^{1}}%
%BeginExpansion
{\displaystyle\int\limits_{0}^{1}}
%EndExpansion
\left\vert h\left(  w\right)  \right\vert \left\vert
\begin{array}
[c]{c}%
\psi^{\prime}\left(  \frac{1+w}{2}\gamma+\frac{1-w}{2}\delta\right)
\end{array}
-\frac{K+k}{2}\right\vert dw\\
+%
%TCIMACRO{\dint \limits_{0}^{1}}%
%BeginExpansion
{\displaystyle\int\limits_{0}^{1}}
%EndExpansion
\left\vert h\left(  w\right)  \right\vert \left\vert
\begin{array}
[c]{c}%
\psi^{\prime}\left(  \frac{1-w}{2}\gamma+\frac{1+w}{2}\delta\right)
\end{array}
-\frac{K+k}{2}\right\vert dw
\end{array}
\right\}  .
\end{align*}

Since $%
\begin{array}
[c]{c}%
\psi^{\prime}%
\end{array}
$ satisfies $-\infty<k\leq%
\begin{array}
[c]{c}%
\psi^{\prime}%
\end{array}
\leq K<+\infty$, we have that%
\[
k-\frac{K+k}{2}\leq\psi^{\prime}-\frac{K+k}{2}\leq K-\frac{K+k}{2},
\]

which implies that
\[
\left\vert \psi^{\prime}-\frac{K+k}{2}\right\vert \leq\frac{K-k}{2}.
\]

Hence,%
\begin{align*}
&  \left\vert
\begin{array}
[c]{c}%
\frac{1}{6}\left[  \psi\left(  \gamma\right)  +4\psi\left(  \frac
{\gamma+\delta}{2}\right)  +\psi\left(  \delta\right)  \right]  -\frac
{2^{\tau-1}}{\left(  \delta-\gamma\right)  ^{\tau}}\frac{\Gamma\left(
\tau+1\right)  }{\Gamma\left(  \tau-m\right)  }\left[
\begin{array}
[c]{c}%
I_{\tau}^{\gamma}\psi\left(  \frac{\gamma+\delta}{2}\right)
\end{array}
+%
\begin{array}
[c]{c}%
^{\delta}I_{\tau}\psi\left(  \frac{\gamma+\delta}{2}\right)
\end{array}
\right] \\
-\frac{K+k}{2}\left(  \delta-\gamma\right)
%TCIMACRO{\dint \limits_{0}^{1}}%
%BeginExpansion
{\displaystyle\int\limits_{0}^{1}}
%EndExpansion
h\left(  w\right)  dw
\end{array}
\right\vert \\
&  \leq\frac{K-k}{2.m!}\frac{\Gamma\left(  \tau+1\right)  }{\Gamma\left(
\tau-m\right)  }\left(  \delta-\gamma\right)
%TCIMACRO{\dint \limits_{0}^{1}}%
%BeginExpansion
{\displaystyle\int\limits_{0}^{1}}
%EndExpansion
\left\vert h\left(  w\right)  \right\vert dw\\
&  \leq\frac{K-k}{2.m!}\frac{\Gamma\left(  \tau+1\right)  }{\Gamma\left(
\tau-m\right)  }\left(  \delta-\gamma\right)  Z_{1}\left(  \tau,m\right)  ,
\end{align*}

where
\[%
%TCIMACRO{\dint \limits_{0}^{1}}%
%BeginExpansion
{\displaystyle\int\limits_{0}^{1}}
%EndExpansion
\left\vert h\left(  w\right)  \right\vert dw=Z_{1}\left(  \tau,m\right)  .
\]

$Z_{1}\left(  \tau,m\right)  $ is defined as in Lemma \ref{3.1}. This ends the proof.
\end{proof}

\begin{remark}
\label{4.2}If we take $\tau=m+1,$ after that if we take $\tau=1$ in Theorem
\ref{4.1}, we obtain the following inequality
\begin{align*}
&  \left\vert
\begin{array}
[c]{c}%
\frac{1}{6}\left[  \psi\left(  \gamma\right)  +4\psi\left(  \frac
{\gamma+\delta}{2}\right)  +\psi\left(  \delta\right)  \right]  -\frac
{1}{\delta-\gamma}%
%TCIMACRO{\dint \limits_{\gamma}^{\delta}}%
%BeginExpansion
{\displaystyle\int\limits_{\gamma}^{\delta}}
%EndExpansion
\psi\left(  \varepsilon\right)
\begin{array}
[c]{c}%
d\varepsilon
\end{array}
\\
-\frac{K+k}{2}\left(  \delta-\gamma\right)
%TCIMACRO{\dint \limits_{0}^{1}}%
%BeginExpansion
{\displaystyle\int\limits_{0}^{1}}
%EndExpansion
h\left(  w\right)
\begin{array}
[c]{c}%
dw
\end{array}
\end{array}
\right\vert \\
&  \leq\frac{5}{72}\left(  K-k\right)  \left(  \delta-\gamma\right)  .
\end{align*}

\end{remark}

Our next aim is an estimation-type result considering the Simpson-like type
conformable fractional integral inequality when $\psi^{\prime}$satisfies a
Lipschitz condition.

\begin{theorem}
\label{4.3}Let $\psi:\left[  \gamma,\delta\right]  \rightarrow%
%TCIMACRO{\U{211d} }%
%BeginExpansion
\mathbb{R}
%EndExpansion
$ be differentiable and continuous mappings on $\left(  a,b\right)  $ and let
$\psi^{\prime}\in L\left[  a,b\right]  .$ Assume that $\psi^{\prime}$
satisfies the Lipschitz condition for some $L>0.$ Then,%
\begin{align}
&  \left\vert
\begin{array}
[c]{c}%
\frac{1}{6}\left[  \psi\left(  \gamma\right)  +4\psi\left(  \frac
{\gamma+\delta}{2}\right)  +\psi\left(  \delta\right)  \right]  -\frac
{2^{\tau-1}}{\left(  \delta-\gamma\right)  ^{\tau}}\frac{\Gamma\left(
\tau+1\right)  }{\Gamma\left(  \tau-m\right)  }\left[
\begin{array}
[c]{c}%
I_{\tau}^{\gamma}\psi\left(  \frac{\gamma+\delta}{2}\right)
\end{array}
+%
\begin{array}
[c]{c}%
^{\delta}I_{\tau}\psi\left(  \frac{\gamma+\delta}{2}\right)
\end{array}
\right] \\
-\left(  \delta-\gamma\right)
\begin{array}
[c]{c}%
\psi^{\prime}\left(  \frac{\gamma+\delta}{2}\right)
\end{array}%
%TCIMACRO{\dint \limits_{0}^{1}}%
%BeginExpansion
{\displaystyle\int\limits_{0}^{1}}
%EndExpansion
h\left(  w\right)  dw
\end{array}
\right\vert \label{4-3}\\
&  \leq L\frac{\left(  \delta-\gamma\right)  ^{2}}{2.m!}\frac{\Gamma\left(
\tau+1\right)  }{\Gamma\left(  \tau-m\right)  }Z_{5}\left(  \tau,m\right)
,\nonumber
\end{align}

\end{theorem}

where
\[
Z_{5}\left(  \tau,m\right)  =%
%TCIMACRO{\dint \limits_{0}^{1}}%
%BeginExpansion
{\displaystyle\int\limits_{0}^{1}}
%EndExpansion
\left\vert h\left(  w\right)  \right\vert wdw
\]
and $h\left(  w\right)  $ is defined as in Theorem \ref{4.1}.

\begin{proof}
From Lemma \ref{3.1}, we have that%
\begin{align*}
&  \frac{1}{6}\left[  \psi\left(  \gamma\right)  +4\psi\left(  \frac
{\gamma+\delta}{2}\right)  +\psi\left(  \delta\right)  \right]  -\frac
{2^{\tau-1}}{\left(  \delta-\gamma\right)  ^{\tau}}\frac{\Gamma\left(
\tau+1\right)  }{\Gamma\left(  \tau-m\right)  }\left[
\begin{array}
[c]{c}%
I_{\tau}^{\gamma}\psi\left(  \frac{\gamma+\delta}{2}\right)
\end{array}
+%
\begin{array}
[c]{c}%
^{\delta}I_{\tau}\psi\left(  \frac{\gamma+\delta}{2}\right)
\end{array}
\right] \\
&  =\frac{\delta-\gamma}{2.m!}\frac{\Gamma\left(  \tau+1\right)  }%
{\Gamma\left(  \tau-m\right)  }\left\{
\begin{array}
[c]{c}%
%TCIMACRO{\dint \limits_{0}^{1}}%
%BeginExpansion
{\displaystyle\int\limits_{0}^{1}}
%EndExpansion
h\left(  w\right)  \left[
\begin{array}
[c]{c}%
\psi^{\prime}\left(  \frac{1+w}{2}\gamma+\frac{1-w}{2}\delta\right)
\end{array}
-%
\begin{array}
[c]{c}%
\psi^{\prime}\left(  \frac{\gamma+\delta}{2}\right)
\end{array}
+%
\begin{array}
[c]{c}%
\psi^{\prime}\left(  \frac{\gamma+\delta}{2}\right)
\end{array}
\right]  dw\\
+%
%TCIMACRO{\dint \limits_{0}^{1}}%
%BeginExpansion
{\displaystyle\int\limits_{0}^{1}}
%EndExpansion
h\left(  w\right)  \left[
\begin{array}
[c]{c}%
\psi^{\prime}\left(  \frac{1-w}{2}\gamma+\frac{1+w}{2}\delta\right)
\end{array}
-%
\begin{array}
[c]{c}%
\psi^{\prime}\left(  \frac{\gamma+\delta}{2}\right)
\end{array}
+%
\begin{array}
[c]{c}%
\psi^{\prime}\left(  \frac{\gamma+\delta}{2}\right)
\end{array}
\right]  dw
\end{array}
\right\} \\
&  =\frac{\delta-\gamma}{2.m!}\frac{\Gamma\left(  \tau+1\right)  }%
{\Gamma\left(  \tau-m\right)  }\left\{
\begin{array}
[c]{c}%
%TCIMACRO{\dint \limits_{0}^{1}}%
%BeginExpansion
{\displaystyle\int\limits_{0}^{1}}
%EndExpansion
h\left(  w\right)  \left[
\begin{array}
[c]{c}%
\psi^{\prime}\left(  \frac{1+w}{2}\gamma+\frac{1-w}{2}\delta\right)
\end{array}
-%
\begin{array}
[c]{c}%
\psi^{\prime}\left(  \frac{\gamma+\delta}{2}\right)
\end{array}
\right]  dw\\
+%
%TCIMACRO{\dint \limits_{0}^{1}}%
%BeginExpansion
{\displaystyle\int\limits_{0}^{1}}
%EndExpansion
h\left(  w\right)  \left[
\begin{array}
[c]{c}%
\psi^{\prime}\left(  \frac{1-w}{2}\gamma+\frac{1+w}{2}\delta\right)
\end{array}
-%
\begin{array}
[c]{c}%
\psi^{\prime}\left(  \frac{\gamma+\delta}{2}\right)
\end{array}
\right]  dw
\end{array}
\right\} \\
&  +\left(  \delta-\gamma\right)
\begin{array}
[c]{c}%
\psi^{\prime}\left(  \frac{\gamma+\delta}{2}\right)
\end{array}%
%TCIMACRO{\dint \limits_{0}^{1}}%
%BeginExpansion
{\displaystyle\int\limits_{0}^{1}}
%EndExpansion
h\left(  w\right)  dw.
\end{align*}

So%
\begin{align*}
&  \left\vert
\begin{array}
[c]{c}%
\frac{1}{6}\left[  \psi\left(  \gamma\right)  +4\psi\left(  \frac
{\gamma+\delta}{2}\right)  +\psi\left(  \delta\right)  \right]  -\frac
{2^{\tau-1}}{\left(  \delta-\gamma\right)  ^{\tau}}\frac{\Gamma\left(
\tau+1\right)  }{\Gamma\left(  \tau-m\right)  }\left[
\begin{array}
[c]{c}%
I_{\tau}^{\gamma}\psi\left(  \frac{\gamma+\delta}{2}\right)
\end{array}
+%
\begin{array}
[c]{c}%
^{\delta}I_{\tau}\psi\left(  \frac{\gamma+\delta}{2}\right)
\end{array}
\right] \\
-\left(  \delta-\gamma\right)
\begin{array}
[c]{c}%
\psi^{\prime}\left(  \frac{\gamma+\delta}{2}\right)
\end{array}%
%TCIMACRO{\dint \limits_{0}^{1}}%
%BeginExpansion
{\displaystyle\int\limits_{0}^{1}}
%EndExpansion
h\left(  w\right)  dw
\end{array}
\right\vert \\
&  \leq\frac{\delta-\gamma}{2.m!}\frac{\Gamma\left(  \tau+1\right)  }%
{\Gamma\left(  \tau-m\right)  }\left\{
\begin{array}
[c]{c}%
%TCIMACRO{\dint \limits_{0}^{1}}%
%BeginExpansion
{\displaystyle\int\limits_{0}^{1}}
%EndExpansion
\left\vert h\left(  w\right)  \right\vert \left\vert
\begin{array}
[c]{c}%
\psi^{\prime}\left(  \frac{1+w}{2}\gamma+\frac{1-w}{2}\delta\right)
\end{array}
-%
\begin{array}
[c]{c}%
\psi^{\prime}\left(  \frac{\gamma+\delta}{2}\right)
\end{array}
\right\vert dw\\
+%
%TCIMACRO{\dint \limits_{0}^{1}}%
%BeginExpansion
{\displaystyle\int\limits_{0}^{1}}
%EndExpansion
\left\vert h\left(  w\right)  \right\vert \left\vert
\begin{array}
[c]{c}%
\psi^{\prime}\left(  \frac{1-w}{2}\gamma+\frac{1+w}{2}\delta\right)
\end{array}
-%
\begin{array}
[c]{c}%
\psi^{\prime}\left(  \frac{\gamma+\delta}{2}\right)
\end{array}
\right\vert dw
\end{array}
\right\}  .
\end{align*}
Since $\psi^{\prime}$satisfies Lipschitz conditions for some $L>0$, we have
that%
\begin{align*}
\left\vert
\begin{array}
[c]{c}%
\psi^{\prime}\left(  \frac{1+w}{2}\gamma+\frac{1-w}{2}\delta\right)
\end{array}
-%
\begin{array}
[c]{c}%
\psi^{\prime}\left(  \frac{\gamma+\delta}{2}\right)
\end{array}
\right\vert  &  \leq L\left\vert \frac{1-w}{2}\gamma+\frac{1+w}{2}\delta
-\frac{\gamma+\delta}{2}\right\vert \\
&  \leq\frac{\delta-\gamma}{2}L\left\vert w\right\vert
\end{align*}
and%
\begin{align*}
\left\vert
\begin{array}
[c]{c}%
\psi^{\prime}\left(  \frac{1-w}{2}\gamma+\frac{1+w}{2}\delta\right)
\end{array}
-%
\begin{array}
[c]{c}%
\psi^{\prime}\left(  \frac{\gamma+\delta}{2}\right)
\end{array}
\right\vert  &  \leq L\left\vert \frac{1-w}{2}\delta+\frac{1+w}{2}\gamma
-\frac{\gamma+\delta}{2}\right\vert \\
&  \leq\frac{\delta-\gamma}{2}L\left\vert w\right\vert .
\end{align*}
Hence,%
\begin{align*}
&  \left\vert
\begin{array}
[c]{c}%
\frac{1}{6}\left[  \psi\left(  \gamma\right)  +4\psi\left(  \frac
{\gamma+\delta}{2}\right)  +\psi\left(  \delta\right)  \right]  -\frac
{2^{\tau-1}}{\left(  \delta-\gamma\right)  ^{\tau}}\frac{\Gamma\left(
\tau+1\right)  }{\Gamma\left(  \tau-m\right)  }\left[
\begin{array}
[c]{c}%
I_{\tau}^{\gamma}\psi\left(  \frac{\gamma+\delta}{2}\right)
\end{array}
+%
\begin{array}
[c]{c}%
^{\delta}I_{\tau}\psi\left(  \frac{\gamma+\delta}{2}\right)
\end{array}
\right] \\
-\left(  \delta-\gamma\right)
\begin{array}
[c]{c}%
\psi^{\prime}\left(  \frac{\gamma+\delta}{2}\right)
\end{array}%
%TCIMACRO{\dint \limits_{0}^{1}}%
%BeginExpansion
{\displaystyle\int\limits_{0}^{1}}
%EndExpansion
h\left(  w\right)  dw
\end{array}
\right\vert \\
&  \leq L\frac{\left(  \delta-\gamma\right)  ^{2}}{2.m!}\frac{\Gamma\left(
\tau+1\right)  }{\Gamma\left(  \tau-m\right)  }%
%TCIMACRO{\dint \limits_{0}^{1}}%
%BeginExpansion
{\displaystyle\int\limits_{0}^{1}}
%EndExpansion
\left\vert h\left(  w\right)  \right\vert wdw\\
&  \leq L\frac{\left(  \delta-\gamma\right)  ^{2}}{2.m!}\frac{\Gamma\left(
\tau+1\right)  }{\Gamma\left(  \tau-m\right)  }Z_{5}\left(  \tau,m\right)  ,
\end{align*}

where
\[%
%TCIMACRO{\dint \limits_{0}^{1}}%
%BeginExpansion
{\displaystyle\int\limits_{0}^{1}}
%EndExpansion
\left\vert h\left(  w\right)  \right\vert wdw=Z_{5}\left(  \tau,m\right)  .
\]
This ends the proof.
\end{proof}

\begin{remark}
\label{4.4}If we take $\tau=m+1,$ after that if we take $\tau=1$ in Theorem
\ref{4.3}, we obtain the following inequality
\begin{align*}
&  \left\vert
\begin{array}
[c]{c}%
\frac{1}{6}\left[  \psi\left(  \gamma\right)  +4\psi\left(  \frac
{\gamma+\delta}{2}\right)  +\psi\left(  \delta\right)  \right]  -\frac
{1}{\delta-\gamma}%
%TCIMACRO{\dint \limits_{\gamma}^{\delta}}%
%BeginExpansion
{\displaystyle\int\limits_{\gamma}^{\delta}}
%EndExpansion
\psi\left(  \varepsilon\right)
\begin{array}
[c]{c}%
d\varepsilon
\end{array}
\\
-\left(  \delta-\gamma\right)
\begin{array}
[c]{c}%
\psi^{\prime}\left(  \frac{\gamma+\delta}{2}\right)
\end{array}%
%TCIMACRO{\dint \limits_{0}^{1}}%
%BeginExpansion
{\displaystyle\int\limits_{0}^{1}}
%EndExpansion
h\left(  w\right)
\begin{array}
[c]{c}%
dw
\end{array}
\end{array}
\right\vert \\
&  \leq L\frac{2\left(  \delta-\gamma\right)  ^{2}}{81},
\end{align*}
with%
\[%
%TCIMACRO{\dint \limits_{0}^{1}}%
%BeginExpansion
{\displaystyle\int\limits_{0}^{1}}
%EndExpansion
\left\vert \frac{1}{3}-\frac{w}{2}\right\vert wdw=\frac{4}{81}.
\]

\end{remark}

\section{Applications}

\subsection{Special Means}

For $0\leq\gamma<\delta,$ we consider the following special means:

\begin{theorem}
\begin{enumerate}
\item[(i)] The arithmetic mean: $A\left(  \gamma,\delta\right)  =\frac
{\gamma+\delta}{2}$,

\item[(ii)] The geometric mean: $G\left(  \gamma,\delta\right)  =\sqrt
{\gamma\delta}$,

\item[(iii)] The logarithmic mean: $L\left(  \gamma,\delta\right)
=\frac{\delta-\gamma}{\ln\delta-\ln\gamma},\gamma\delta\neq0,$

\item[(iv)] The logarithmic mean: $L_{s}\left(  \gamma,\delta\right)  =\left(
\frac{\delta^{s+1}-\gamma^{s+1}}{\left(  s+1\right)  \left(  \gamma
-\delta\right)  }\right)  ^{\frac{1}{s}},s\in%
%TCIMACRO{\U{2124} }%
%BeginExpansion
\mathbb{Z}
%EndExpansion
\backslash\left\{  0,1\right\}  .$
\end{enumerate}
\end{theorem}

Next, using the main results obtained in Section 3, we give some applications
to special means of real numbers.

\begin{proposition}
\label{5.1}Let $0<\gamma<\delta,s\in%
%TCIMACRO{\U{2115} }%
%BeginExpansion
\mathbb{N}
%EndExpansion
$ $.$ Then%
\begin{align}
&  \left\vert \frac{1}{3}A\left(  \gamma^{s},\delta^{s}\right)  +\frac{2}%
{3}A^{s}\left(  \gamma,\delta\right)  -L_{s}^{s}\left(  \gamma,\delta\right)
\right\vert \label{5-1}\\
&  \leq\frac{5\left(  \delta-\gamma\right)  }{72}s\left[  \gamma^{s-1}%
+\delta^{s-1}\right]  .\nonumber
\end{align}

\end{proposition}

\begin{proof}
The proof is obvious from Remark \ref{3.4} applied $\psi\left(  \varepsilon
\right)  =\varepsilon^{s}.$
\end{proof}

\begin{remark}
\label{5.2a}If we take $s=1$ in \ref{5-1}, we obtain the following inequality
\begin{equation}
\left\vert A\left(  \gamma,\delta\right)  -L\left(  \gamma,\delta\right)
\right\vert \leq\frac{5\left(  \delta-\gamma\right)  }{36}. \label{5-2a}%
\end{equation}

\end{remark}

See also (\cite{ADD09}, Page 13).

\begin{proposition}
\label{5.3}Let $0<\gamma<\delta,s\in%
%TCIMACRO{\U{2115} }%
%BeginExpansion
\mathbb{N}
%EndExpansion
$ $.$ Then%
\begin{align}
&  \left\vert \frac{1}{3}A\left(  \gamma^{s},\delta^{s}\right)  +\frac{2}%
{3}A^{s}\left(  \gamma,\delta\right)  -L_{s}^{s}\left(  \gamma,\delta\right)
\right\vert \label{5-3}\\
&  \leq\frac{\delta-\gamma}{2}\left(  \frac{2^{p+1}+1}{3\left(  p+1\right)
}\right)  ^{\frac{1}{p}}\left(  \frac{1}{4}\right)  ^{\frac{1}{q}}\left[
\begin{array}
[c]{c}%
\left(  \left(  s\gamma^{s-1}\right)  ^{q}+3\left(  s\delta^{s-1}\right)
^{q}\right)  ^{\frac{1}{q}}\\
+\left(  3\left(  s\gamma^{s-1}\right)  ^{q}+\left(  s\delta^{s-1}\right)
^{q}\right)  ^{\frac{1}{q}}%
\end{array}
\right]  .\nonumber
\end{align}

\end{proposition}

\begin{proof}
The proof is obvious from Remark \ref{3.6a} applied $\psi\left(
\varepsilon\right)  =\varepsilon^{s}.$
\end{proof}

\begin{remark}
\label{5.3a}If we take $s=1$ in \ref{5-3}, we obtain the following inequality%
\begin{equation}
\left\vert A\left(  \gamma,\delta\right)  -L\left(  \gamma,\delta\right)
\right\vert \leq\left(  \delta-\gamma\right)  \left(  \frac{2^{p+1}%
+1}{3\left(  p+1\right)  }\right)  ^{\frac{1}{p}}. \label{5-3a}%
\end{equation}

\end{remark}

\begin{proposition}
\label{5.4}Let $0<\gamma<\delta,s\in%
%TCIMACRO{\U{2115} }%
%BeginExpansion
\mathbb{N}
%EndExpansion
$ and $0<q<1.$ Then%
\begin{align}
&  \left\vert \frac{1}{3}A\left(  \gamma^{s},\delta^{s}\right)  +\frac{2}%
{3}A^{s}\left(  \gamma,\delta\right)  -L_{s}^{s}\left(  \gamma,\delta\right)
\right\vert \label{5-4}\\
&  \leq\frac{\delta-\gamma}{2}\left(  \frac{5}{36}\right)  ^{1-\frac{1}{q}%
}\left(  \frac{1}{648}\right)  ^{\frac{1}{q}}\left[
\begin{array}
[c]{c}%
\left(  61\left(  s\gamma^{s-1}\right)  ^{q}+29\left(  s\delta^{s-1}\right)
^{q}\right)  ^{\frac{1}{q}}\\
+\left(  29\left(  s\gamma^{s-1}\right)  ^{q}+61\left(  s\delta^{s-1}\right)
^{q}\right)  ^{\frac{1}{q}}%
\end{array}
\right]  .\nonumber
\end{align}

\end{proposition}

\begin{proof}
The proof is obvious from Remark \ref{3.8} applied $\psi\left(  \varepsilon
\right)  =\varepsilon^{s}.$
\end{proof}

\begin{remark}
\label{5.5}If we take $s=1$ in \ref{5-4}, we obtain the following inequality%
\begin{equation}
\left\vert A\left(  \gamma,\delta\right)  -L\left(  \gamma,\delta\right)
\right\vert \leq\frac{5\left(  \delta-\gamma\right)  }{36}. \label{5-5}%
\end{equation}

\end{remark}

See also (\cite{ADD09}, Page 13).

\begin{proposition}
\label{5.6}Let $0<\gamma<\delta,s\in%
%TCIMACRO{\U{2115} }%
%BeginExpansion
\mathbb{N}
%EndExpansion
.$ Then%
\begin{align}
&  \left\vert \frac{1}{3}A\left(  \gamma^{s},\delta^{s}\right)  +\frac{2}%
{3}A^{s}\left(  \gamma,\delta\right)  -L_{s}^{s}\left(  \gamma,\delta\right)
\right\vert \label{5-6}\\
&  \leq\frac{\delta-\gamma}{12}\left(  \frac{2^{p+1}+1}{3\left(  p+1\right)
}\right)  ^{\frac{1}{p}}\left(  \frac{1}{4}\right)  ^{\frac{1}{q}}\left[
\begin{array}
[c]{c}%
\left(  \left(  s\gamma^{s-1}\right)  ^{q}+3\left(  s\left(  \frac
{\gamma+\delta}{2}\right)  ^{s-1}\right)  ^{q}\right)  ^{\frac{1}{q}}\\
+\left(  3\left(  s\left(  \frac{\gamma+\delta}{2}\right)  ^{s-1}\right)
^{q}+\left(  s\delta^{s-1}\right)  ^{q}\right)  ^{\frac{1}{q}}%
\end{array}
\right]  .\nonumber
\end{align}

\end{proposition}

\begin{remark}
\label{5.7}If we take $s=1$ in \ref{5-6}, we obtain the following inequality%
\begin{equation}
\left\vert A\left(  \gamma,\delta\right)  -L\left(  \gamma,\delta\right)
\right\vert \leq\frac{\delta-\gamma}{6}\left(  \frac{2^{p+1}+1}{3\left(
p+1\right)  }\right)  ^{\frac{1}{p}}. \label{5-7}%
\end{equation}

\end{remark}

\begin{proposition}
\label{5.8} Let $0<\gamma<\delta.$ Then%
\begin{align}
&  \left\vert \frac{1}{3}A\left(  \alpha,\beta\right)  +\frac{2}{3}G\left(
\alpha,\beta\right)  -L\left(  \alpha,\beta\right)  \right\vert \label{5-8}\\
&  \leq\left(  \ln\beta-\ln\alpha\right)  \frac{5}{36}A\left(  \alpha
,\beta\right)  .\nonumber
\end{align}

\end{proposition}

\begin{proof}
The proof is obvious from the Remark \ref{3.4}, applied $f\left(
\varepsilon\right)  =e^{\varepsilon},\varepsilon>0$ and $\alpha=e^{\gamma
},\beta=e^{\delta}.$
\end{proof}

\subsection{Inequalities for some special functions.}

\subsubsection{Modified Bessel function.}

Recall the first kind modified Bessel function $I_{\rho}$ , which has the
series representation (\cite{W44}, p.77)
\[
I_{\rho}\left(  \varepsilon\right)  =\sum_{n\geq0}\frac{\left(  \frac
{\varepsilon}{2}\right)  ^{\rho+2n}}{n!\Gamma\left(  \rho+n+1\right)  },
\]
where $\varepsilon\in%
%TCIMACRO{\U{211d} }%
%BeginExpansion
\mathbb{R}
%EndExpansion
$ and $\rho>-1,$ while the second kind modified Bessel function $K\rho$
(\cite{W44}, p.78) is usually defined as%
\[
K\rho\left(  \varepsilon\right)  =\frac{\pi}{2}\frac{I_{-\rho}\left(
\varepsilon\right)  -I_{\rho}\left(  \varepsilon\right)  }{\sin\rho\pi}.
\]

Here , we consider the function $\wp_{\rho}:%
%TCIMACRO{\U{211d} }%
%BeginExpansion
\mathbb{R}
%EndExpansion
\rightarrow\left[  1,\infty\right)  $ defined by%
\[
\wp_{\rho}\left(  \varepsilon\right)  =2^{\rho}\Gamma\left(  \rho+1\right)
\varepsilon^{-\rho}I_{\rho}\left(  \varepsilon\right)  ,
\]

where $\Gamma$ is the Gamma function.

\begin{proposition}
\label{5.9} Let $\rho>-1$ and $0<\gamma<\delta.$ Then%
\begin{align}
&  \left\vert \frac{1}{6}\left[  \wp_{\rho}\left(  \gamma\right)  +4\wp_{\rho
}\left(  \frac{\gamma+\delta}{2}\right)  +\wp_{\rho}\left(  \delta\right)
\right]  -\frac{1}{\delta-\gamma}%
%TCIMACRO{\dint \limits_{\gamma}^{\delta}}%
%BeginExpansion
{\displaystyle\int\limits_{\gamma}^{\delta}}
%EndExpansion
\wp_{\rho}\left(  \varepsilon\right)  d\varepsilon\right\vert \label{5-9}\\
&  \leq\frac{\delta-\gamma}{2\left(  \left(  \rho+1\right)  \right)  }\left(
\frac{5}{36}\right)  ^{1-\frac{1}{q}}\left(  \frac{1}{648}\right)  ^{\frac
{1}{q}}\left(
\begin{array}
[c]{c}%
\left(  61\gamma\left\vert \wp_{\rho+1}\left(  \gamma\right)  \right\vert
^{q}+29\delta\left\vert \wp_{\rho+1}\left(  \delta\right)  \right\vert
^{q}\right)  ^{\frac{1}{q}}\\
+\left(  29\gamma\left\vert \wp_{\rho+1}\left(  \gamma\right)  \right\vert
^{q}+61\delta\left\vert \wp_{\rho+1}\left(  \delta\right)  \right\vert
^{q}\right)  ^{\frac{1}{q}}%
\end{array}
\right)  .\nonumber
\end{align}

\end{proposition}

Specially, if $\rho=\frac{-1}{2},$ then%
\begin{align}
&  \left\vert \frac{1}{6}\left[  \cosh\left(  \gamma\right)  +4\cosh\left(
\frac{\gamma+\delta}{2}\right)  +\cosh\left(  \delta\right)  \right]
-\frac{1}{\delta-\gamma}%
%TCIMACRO{\dint \limits_{\gamma}^{\delta}}%
%BeginExpansion
{\displaystyle\int\limits_{\gamma}^{\delta}}
%EndExpansion
\cosh\left(  \varepsilon\right)  d\varepsilon\right\vert \label{5-10}\\
&  \leq\frac{\delta-\gamma}{2\left(  \left(  \rho+1\right)  \right)  }\left(
\frac{5}{36}\right)  ^{1-\frac{1}{q}}\left(  \frac{1}{648}\right)  ^{\frac
{1}{q}}\left(
\begin{array}
[c]{c}%
\left(  61\gamma\left\vert \frac{\sinh\left(  \gamma\right)  }{\gamma
}\right\vert ^{q}+29\delta\left\vert \frac{\sinh\left(  \delta\right)
}{\delta}\right\vert ^{q}\right)  ^{\frac{1}{q}}\\
+\left(  29\gamma\left\vert \frac{\sinh\left(  \gamma\right)  }{\gamma
}\right\vert ^{q}+61\delta\left\vert \frac{\sinh\left(  \delta\right)
}{\delta}\right\vert ^{q}\right)  ^{\frac{1}{q}}%
\end{array}
\right)  .\nonumber
\end{align}

\begin{proof}
Apply inequality \ref{3-8} to the mapping $f\left(  \varepsilon\right)
=\wp_{\rho}\left(  \varepsilon\right)  ,\varepsilon>0$ and $\wp_{\rho}%
^{\prime}\left(  \varepsilon\right)  =\frac{\varepsilon}{\rho+1}\wp_{\rho
+1}\left(  \varepsilon\right)  .$ Now taking into account the relations
$\wp_{\frac{-1}{2}}\left(  \varepsilon\right)  =\cosh\left(  \varepsilon
\right)  $ and $\wp_{\frac{1}{2}}\left(  \varepsilon\right)  =\frac
{\sinh\left(  \varepsilon\right)  }{\varepsilon},$ we have that inequality
\ref{5-9} is reduced to inequality \ref{5-10}.
\end{proof}

\subsubsection{$q-$digamma function}

Let $0<q<1.$ The $q-$digamma function $\Psi_{q}$ is the $q-$analogue of the
$\Psi$ or digamma function $\Psi$ defined by%
\begin{align*}
\Psi_{q}\left(  \eta\right)   &  =-\ln\left(  1-q\right)  +\ln q%
%TCIMACRO{\dsum \limits_{u=0}^{\infty}}%
%BeginExpansion
{\displaystyle\sum\limits_{u=0}^{\infty}}
%EndExpansion
\frac{%
\begin{array}
[c]{c}%
q
\end{array}
^{u+\eta}}{1-%
\begin{array}
[c]{c}%
q
\end{array}
^{u+\eta}}\\
&  =-\ln\left(  1-q\right)  +\ln q%
%TCIMACRO{\dsum \limits_{u=0}^{\infty}}%
%BeginExpansion
{\displaystyle\sum\limits_{u=0}^{\infty}}
%EndExpansion
\frac{%
\begin{array}
[c]{c}%
q
\end{array}
^{u\eta}}{1-%
\begin{array}
[c]{c}%
q
\end{array}
^{u}}.
\end{align*}
For $q>1$ and $\eta>0,$ the $q-$digamma function $\Psi_{q}$ is defined by%
\begin{align*}
\Psi_{q}\left(  \eta\right)   &  =-\ln\left(  q-1\right)  +\ln q\left(
\eta-\frac{1}{2}-%
%TCIMACRO{\dsum \limits_{u=0}^{\infty}}%
%BeginExpansion
{\displaystyle\sum\limits_{u=0}^{\infty}}
%EndExpansion
\frac{%
\begin{array}
[c]{c}%
q
\end{array}
^{-\left(  u+\eta\right)  }}{1-%
\begin{array}
[c]{c}%
q
\end{array}
-^{\left(  u+\eta\right)  }}\right) \\
&  =-\ln\left(  q-1\right)  +\ln q\left(  \eta-\frac{1}{2}-%
%TCIMACRO{\dsum \limits_{u=0}^{\infty}}%
%BeginExpansion
{\displaystyle\sum\limits_{u=0}^{\infty}}
%EndExpansion
\frac{%
\begin{array}
[c]{c}%
q
\end{array}
^{-u\eta}}{1-%
\begin{array}
[c]{c}%
q
\end{array}
-^{u\eta}}\right)  .
\end{align*}

\begin{proposition}
\label{5.11}Let $0<\gamma<\delta,q>1$ and $\frac{1}{p}+\frac{1}{q}=1.$ Then
\begin{align}
&  \left\vert \frac{1}{6}\left[  \Psi_{q}\left(  \gamma\right)  +4\Psi
_{q}\left(  \frac{\gamma+\delta}{2}\right)  +\Psi_{q}\left(  \delta\right)
\right]  -\frac{1}{\delta-\gamma}%
%TCIMACRO{\dint \limits_{\gamma}^{\delta}}%
%BeginExpansion
{\displaystyle\int\limits_{\gamma}^{\delta}}
%EndExpansion
\Psi_{q}\left(  \varepsilon\right)  d\varepsilon\right\vert \label{5-11}\\
&  \leq\frac{\delta-\gamma}{12}\left(  \frac{2^{p+1}+1}{3\left(  p+1\right)
}\right)  ^{\frac{1}{p}}\left(  \frac{1}{4}\right)  ^{\frac{1}{q}}\left[
\begin{array}
[c]{c}%
\left(  \left\vert \Psi_{q}^{\prime}\left(  \gamma\right)  \right\vert
^{q}+3\left\vert \Psi_{q}^{\prime}\left(  \frac{\gamma+\delta}{2}\right)
\right\vert ^{q}\right)  ^{\frac{1}{q}}\\
\left(  3\left\vert \Psi_{q}^{\prime}\left(  \frac{\gamma+\delta}{2}\right)
\right\vert ^{q}+\left\vert \Psi_{q}^{\prime}\left(  \delta\right)
\right\vert ^{q}\right)  ^{\frac{1}{q}}%
\end{array}
\right]  .\nonumber
\end{align}

\end{proposition}

\begin{proof}
The assertion can be obtained immediately by using Remark \ref{3.10} to
$\psi\left(  \varepsilon\right)  =\Psi_{q}\left(  \varepsilon\right)  ,$ and
$\varepsilon>0,$ since $\psi^{\prime}\left(  x\right)  =\Psi_{q}^{\prime
}\left(  \varepsilon\right)  $ is convex on $\left(  0,+\infty\right)  .$
\end{proof}

\section{Conclusion}

In this paper, using a new identity of Simpson-like type for conformable
fractional integral, we obtained some new Simpson type conformable fractional
integral inequalities. Furthermore, some interesting applications were
examined, for example, we applied the investigated results to special means of
real numbers and two special functions named modified Bessel function and
$q-$digamma function. So, this paper is a detailed examination of the
Simpson-like type conformable fractional integral inequalities.

\end{document}